\documentclass[a4 paper]{article}

\usepackage[utf8]{inputenc}
\usepackage[french]{babel} 
\usepackage[T1]{fontenc} 

\usepackage{amsmath}
\usepackage{amssymb}
\usepackage{color}
\usepackage{amscd}
\usepackage{amsfonts}
\usepackage{graphicx}

\usepackage{fullpage}

\usepackage{setspace}  

\newcommand{\vx}{\vspace{0,5cm}}

\newcommand{\zp}{\mathbb{Z}_p}
\newcommand{\zd}{\mathbb{Z}_2}
\newcommand{\z}{\mathbb{Z}}

\newcommand{\ft}{F^\times}
\newcommand{\kt}{K^\times}
\newcommand{\gal}{\text{Gal}}
\newcommand{\qp}{\mathbb{Q}_p}

\begin{document}

 \begin{spacing}{1.7}

	\begin{center}	
		{\Large {\bf Galois module structure of cyclic extensions of local fields of characteristic zero}}

		\vspace{0.1cm}
		
		{ {\textsc{ Sébastien Bosca }}}
	\end{center}

\vspace{1cm}


\vspace{0,5cm}
 
 {\bf Abstract:} Let $K$ be a non archimedean local field of characteristic zero. $F/K$ being a cyclic extension of degree $p^n$, we determine the $\mathbb{Z}_p[G]$-module $F^\times$ up to isomorphism.
 
\vspace{0,5cm} 


\vspace{1cm}

\end{spacing}

{\bf \large I Introduction} \vspace{1cm}

Class Field Theory describes abelian extensions of number fields but raises at least as many questions as it answers. We just still don't know precisely when the fundamental unit of $\mathbb{Q}[\sqrt{pq}]$ has a norm equals to 1 or -1, see [Schol] also discussed by Lemmermeyer in [Lem]. Basic facts about how $K^\times$ maps into its idele group are unknown (see Gross, Leopold, Iwasawa conjectures). 

Iwasawa Theory usually studies $G$-modules up to pseudo-isomorphisms, not isomorphisms.
On the local side, obviously easier than the global one, the character of $\qp\underset{\mathbb{Z}_p}{\otimes}F^\times$ for any Galois extension $F/K$ is known, due to the normal basis theorem and the log application. But this character is not enough to determine completely the $G$-module $F^\times$. This leaves room for finite differences. Moving along $K_n/K$ in a $\mathbb{Z}_p$-extension, these finite differences in the idele group of $K_n$ may increase with $n$ and lead to an unexpected $\lambda$ invariant for classical $\Lambda$-modules of $K_\infty$.

So, what has been done about the structure of $F^\times$ for a given local cyclic extension $F/K$ ? When the order of the Galois group of $F/K$ is $p^n$ and when $K$ has no roots of unity with order a power of $p$, a result from Yakovlev issued in 1970 allows to determine the $G$-module $F^\times$, see [Yak]. More recently, Sharifi studied local unit groups using filtrations, see [Sha].

In this context, the main motivation of this paper was to determine whether basic parameters are enough or not to describe completely the $G$-module $F^\times$ in the cyclic case, the simplest one. The answer given is that for any cyclic extension of local fields $F/K$ of degree $p^n$, the dimension of $K$ over $\qp$, the fact that the residual characteritic is $p$ or not, the number of roots of unity in $K$, in $F$, the way $\gal(F/K)$ acts on them and how much are norm in $F/K$ is enough to know $F^\times$ up to isomorphism.

\vspace{1cm}

{\bf \large II Basic tools and notations} \vspace{1cm}

Throughout this paper, $K$ denotes a local $p$-adic field, and $F$ a cyclic extension of $K$ of degree $p^n$ for some integer $n> 0$, $G$ being its Galois group and $\sigma$ a generator of $G$.

We only deal with $p$-parts of $K^\times$, $F^\times$ and so on, using $p$-completion as described at the beginning of III.

$\mu_K$ denotes the group of roots of unity in $K$, of which the order is a power of $p$.

In any $A$-module $M$, $\langle x,y\rangle$ denotes the submodule generated by $x$ and $y$; when $A$ may not be clear we specify it, notably $\langle x,y\rangle _{\zp}$ and $\langle x,y\rangle _{\zp[G]}$ respectively denote the $\zp$-module and the $\zp[G]$-module generated by $x$ and $y$.

\newpage

Moreover, we will use concepts and toolds coming from the following areas:\vspace{0,5cm}

\begin{itemize}
\item[\hspace{1cm}(1.1)] Kummer theory: radical of $F/K$ if $K$ contains roots of unity, see [Mil], appendix p226 or [CaFr], chapter III pp 85-93. The fact that for a given cyclic extension $F/K$, $L/F$ being a Kummer extension, $G=\gal(F/K)$ acts $\omega$-isotypically on $Rad(L/F)$ iff $L\subset K^{ab}$ ($K^{ab}$ being the maximal abelian extension of $K$ and $\omega$ being the character of the action of $G$ on $\mu_F$) will be used in III, Case 5.   \vspace{0,5cm}

\item[\hspace{1cm}(1.2)] Class field theory: if $F/K$ is a cyclic extension of local fields,
$$\frac{\kt}{N_{F/K}(\ft)}=\gal(F/K)\ ,$$
and the compatibility between restriction of Galois group and the norm application (see [CaFr] or [Mil]) stating that the following diagram is commutative if $K\subset L$ are 2 local fields, $\kt$ and $\ft$ the completions of their multiplicative groups:
$$\begin{matrix}
L^\times & \widetilde{\longrightarrow} &\gal(L^{ab}/L)\\
 & & \\
| &  &| \\
N_{L/K} & & Res\\
\big\downarrow& &\big\downarrow\\
 & & \\ 
\kt & \widetilde{\longrightarrow} &\gal(K^{ab}/K)
\end{matrix}$$

\vspace{0,5cm}

\item[\hspace{1cm} (1.3)]  $N(L)=N(M)$ lemma:\vspace{0,5cm}

Suppose that $M$ is a $\zp[G]$-module of finite type and that $$H^1(G,M)=\{0\}\ ;$$
and suppose that $L\subset M$ satisfies $$N(L)=N(M)$$

where $N=N_{F/K}=1+\sigma+\sigma^2+...+\sigma^{p^n-1}$ is the norm application, then, one has:
$$L=M$$
\vspace{0,5cm}

Proof: $H^1(G,M)$ means that $M^{1-\sigma}$ is the kernel of the norm application ($\sigma$ being a generator of $G$), so that $N(L)=N(M)$ means that $M=LM^{1-\sigma}$; then, $Q=M/L$ satisfies $Q=Q^{1-\sigma}$ and you can apply Nakayama's lemma in $Q$ as a $\zp[G]$-module, $\zp[G]$ being a local ring with maximal ideal $I=(1-\sigma,p)$.\vspace{0,5cm}

\item[\hspace{1cm}(1.4)] Cohomology: definitions  of $H^0(G,M)$, $H^1(G,M)$ where $G$ is a cyclic group and $M$ a $\zp[G]$-module, basic properties of the Herbrand quotient: $$q(G,M)=\frac{|H^0(G,M)|}{|H^1(G,M)|}$$

see [Lang], chapter IX or [Mil] p 90 and, for a use case [Bos] p 9.\vspace{0,5cm}

\item[\hspace{1cm}(1.5)] Characters of finite groups: regular character, properties related to direct sums and quotients (see [Ser]), and the following result:

 \begin{center}
 {\bf Property:}\vspace{0,5cm}

Let $\chi$ be a $\qp$-irreducible character of $G=\mathbb{Z}/p^n\mathbb{Z}$; then its Herbrand quotient is: $$q(\chi)=\left\{\begin{matrix}
p^n&\ if &\chi = 1\\
\frac{1}{p} &\ if &\chi \neq 1
\end{matrix}\right.$$
 \end{center}

where $1$ means the character of $\zp$ in which $G$ acts trivially.

 To prove it, note that any monogene $\qp[G]$-module is a quotient of $\frac{\qp[X]}{X^{p^n}-1}=\prod_{k=0}^n \frac{\qp[X]}{P_k(X)}=\prod V_k$ where the polynomials $P_k(X)$ are $P_k(X)=\frac{X^{p^k}-1}{X^{p^{k-1}}-1}$ for $k> 0$ and $P_0(X)=X-1$; each $V_k$ is irreducible and an easy computation gives $q(V_0)=p^n$, and $q(V_k)=\frac{1}{p}$ if $k> 0$.

\vspace{0,5cm}

\item[\hspace{1cm}(1.6)] $N\in\zp[G]$ or $N_{F/K}$ denotes the norm application:
$$N=1+\sigma+\sigma^2+....+\sigma^{p^n-1}=\sum_{k\in\mathbb{Z}/p^n\mathbb{Z}}\sigma^k$$

\item[\hspace{1cm}(1.7)] $A\in\zp[G]$ denotes the following:
$$A=0+1\sigma+2\sigma^2+....+(p^n-1)\sigma^{p^n-1}=\sum_{\ \ \ 0\le k\le {p^n-1}}k\sigma^k$$
and has the following property in $\zp[G]$:
$$A(1-\sigma)=N-p^n$$

Proof: that's an Abel transform: $$A(1-\sigma)=\sigma+2\sigma^2+...+(p^n-1)\sigma^{p^n-1}-\sigma^2-2\sigma^3-...-(p^n-2)\sigma^{p^n-1}-(p^n-1)\sigma^{p^n}$$
$$=(1-0)\sigma+(2-1)\sigma^2+...+(p^n-1-(p^n-2))\sigma^{p^n-1}-(p^n-1).1$$
$$=\sigma+\sigma^2+...+\sigma^{p^n-1}-p^n+1$$
$$=N-p^n\ .$$

\vspace{0,5cm}

\item[\hspace{1cm}(1.8)] A basic result about $\zp$-modules:\vspace{0,5cm}

 Let $M$ be a $\zp$-module with no torsion of finite type. Then, for any submodule $L\subset M$, the following conditions are equivalent:\vspace{0,3cm}
 \begin{enumerate}
  	\item$L$ is a direct factor in $M$: $M=L\oplus S$ for some $S$
 	\item \hspace{1cm} $L\cap M^p=L^p$\vx
 	
 	More precisely, in order to prove that $L$ is a direct factor in $M$, for a given family $l_1,...,l_k$ of elements generating $L$, if we prove
 	
 	$$\prod_{i=1}^{k}l_i^{\alpha_i}\in M^p\ \iff \ \forall i\in [\![1;k]\!]:\ p\mid\alpha_i$$
 	
 	we prove both that $L$ is a direct factor in $M$ and that the $(l_i)_{1\le i\le k}$ are linearly independent in $L/L^p$.
 	
 \end{enumerate}

\end{itemize}\vspace{2cm}

{\bf \large III Computations in cyclic extensions} \vspace{1cm}

Let $F/K$ be an extension of local fields of characteristic $0$. We note $G=\text{Gal}(F/K)$, generated by $\sigma$. As above, $|G|=p^n$ is the number of elements in $G$. Throughout, we replace the modules $K^\times$, $F^\times$ and so on tacitly by their
p-completions. The p-completion of a finite module is simply its p-part, so we will also say p-part for p-completion in general. 

$$\overline{K^\times}=\underset{\underset{n\in\mathbb{N}}{\longleftarrow}}{\ \lim\ }\ \ \frac{K^\times}{K^{\times p^n}}$$

and the same inverse limit $\overline{F^\times}$ holds for $F^\times$. What is the difference between $K^\times$ and $\overline{K^\times}$ ? Essentially, we then replace the subgroup $\pi^\mathbb{Z}$ generated by a uniformizing parameter with, say, $\pi^{\zp}$; and we remove all roots of unity having an order prime to $p$.

However these notations are not of high interest here and we'll omit them, going on using $K^\times$ and $F^\times$ instead of $\overline{K^\times}$ and $\overline{F^\times}$.

When calculating the character of a $\zp[G]$-module $M$, we mean the one of $\mathbb{Q}_p\underset{\zp}{\otimes}M$. 1 denotes the character of $\zp$ on which $G$ acts trivially, and $\chi_{reg}$ the one of $\zp[G]$.

Sometimes we'll need to consider subfields of $F$ containing $K$, and then $K_i$ denotes the unique field $K\subset K_i\subset F$ with $[K_i/K]=p^i$:

$$K=K_0\subset K_1\subset K_2\subset ... \subset K_n=F$$

We will gradually use the following parameters:
\begin{enumerate}
	\item Denoting $\mu_K$ the set of roots of unity in $K$ (with order a power of $p$), $a$ is the positive integer such that $|\mu_K|=p^a$; also, $\mu_{p^k}$ is the group of roots of unity of order $p^i$, $0\le i\le k$ in an algebraic closure of $K$ and $\mu_\infty=\underset{k\in\mathbb{N}}{\cup}\mu_K$.
	\item In the cyclic extension $F/K$ of order $p^n$, $b$ is the positive integer such that the order of $\mu_K$ in $\kt/N(\ft)$ is $p^b$. 
	Note that one has both $b\le a$ and $b\le n$ because $p^b$ is the order of an element of order $p^a$ in a group of order $p^n$.

	\item $m$ is the positive integer such that $K_m=K[\mu_F]$ is the maximal cyclotomic extension of $K$ contained in $F$. If $p\neq 2$ or $a\ge 2$ it means that $K_m=F\cap K[\mu_\infty]$; however, if $p=2$ and $i\notin K$, this equality may fail when $K[\mu_\infty]/K$ is not procyclic, as explained below.
\end{enumerate}

\vx

In order to prove the main result, which is theorem IV.5, we'll have to cover many cases. We begin with the simpliest ones, while the 5th is the most important. The reader will have to choose which ones are interesting for him.

\vspace{1cm}
{\bf \large Case 1: $\mu_F=\mu_K=\{1\}$} \vspace{1cm}

This is the most simple case, which can be seen as a consequence of results proved by Yakovlev in [Yak]. The $\log$ application is an isomorphism between an open subgroup of the unit group $\mathcal{U}_F$ of the local field $F$, and an open subgroup of the group $\mathcal{O}_F$ of integers of $F$. The normal basis theorem states that the character of $F/K$ is the regular one, so 

$$\chi_{\mathcal{O}_F}=\chi_\text{reg}$$

as an $\mathcal{O}_K[G]$-module and then,
$$\chi_{\mathcal{O}_F}=d\chi_\text{reg}$$
as a $\zp[G]$-module, where $d=[K:\mathbb{Q}_p]$. Hence,
$$\chi_{\mathcal{U_F}}=d\chi_\text{reg}\ ,$$
so that $$\chi_{{F^\times}}=d\chi_\text{reg}+1\ .$$

\vspace{0,5cm}

Class Field Theory states that $$N_{F/K}(F^\times)$$

is an open subgroup of $K^\times$ of index $|G|=p^n$.

So, we choose $x\in K^\times$ of which the order in $K^\times/N(F^\times)$ is $p^n$, and we choose a family $v_1,...,v_d$ in $N(F^\times)$ such that $v_1,...,v_d$ is a basis of $N(F^\times)$ in the $\mathbb{F}_p$-vector space $K^\times/K^{\times p}$ of dimension $d+1$ (because $K$ has no roots of unity).
$$x\notin N(F^\times)K^{\times p}$$
and then, 
$$\kt/\langle v_1,...,v_d\rangle $$ is cyclic and generated by $x$, so that 
$$\kt/\langle v_1,...,v_d, x^{p^n}\rangle $$
is cyclic with order less than or equal to $p^n$; but now 
$$\langle v_1,...,v_d, x^{p^n}\rangle \subset N(\ft)$$ and finally, one has

$$N(\ft)=\langle v_1,...,v_d,x^{p^n}\rangle \ .$$

Now we can find a family $u_1,...,u_d$ in $F^\times$ such that $N(u_i)=v_i$ for $i=1,...,d$, and one has $N(x)=x^{p^n}$. The "$N(L)=N(M)$" lemma concludes that

$$\ft=\langle u_1,...,u_d,x\rangle \ .$$

Since we know that $x$ is invariant, $\ft$ is then a quotient of $$\zp[G]^d\times\zp\ ;$$
but its character is precisely $d\chi_{reg}+1$, so we are sure that the $(u_i)_{1\le i\le d}$ generate a free $\zp[G]$-module of rank $d$ and that 

$$\ft\simeq\zp[G]^d\times\zp$$

\vspace{1cm}
{\bf \large Case 2: $\mu_F=\mu_K=\mu_{p^a}\subset N(\ft)$} 
\vspace{1cm}

For this more complex case, we'll use the 2 following lemma:\vspace{0,5cm}

{\bf lemma (3.2.1):} computing $N(y)$ when $z^{p^a}=y^{1-\sigma}, N(z)=\xi$\vspace{0,5cm}

Let $F/K$ be a finite cyclic extension of fields (local or not) of degree $p^n$, with $\mu_{p^a}\subset K$. Then, the relation $$z^{p^a}=y^{1-\sigma}$$ for $z$ satisfying $N(z)=\xi$ where $\xi$ is a primitive root of unity of order $p^a$, implies that, for $c=\min (n,a)$, one has:   

$$r=z^{p^{a-c}A}y^{p^{n-c}}=\left\{\begin{matrix}
	z^{p^{a-n}A}y& \text{ if }&n\le a
	\\z^Ay^{p^{n-a}}& \text{ if }&n\ge a
\end{matrix}\right.\ ,$$

is such that $$r^{(1-\sigma)}=\xi^{p^{a-c}}$$
so that $$r_0=r^{p^c}\in K$$
is the radical of the unique extension of degree $p^c$ contained in $F$. Moreover, one has:
$$N(y)=r_0$$

\vspace{0,5cm}

Proof:\vspace{0,5cm}

We use the formula exposed in (1.7), $A(1-\sigma)=N-p^n$:\vspace{0,5cm}

$$r=z^{p^{a-c}A}y^{p^{n-c}}$$
$$r^{(1-\sigma)}=z^{p^{a-c}A(1-\sigma)}y^{p^{n-c}(1-\sigma)}$$
$$r^{(1-\sigma)}=\frac{N(z)^{p^{a-c}}}{z^{p^{n+a-c}}}z^{p^{a+n-c}}$$
$$r^{(1-\sigma)}=N(z)^{p^{a-c}}=\xi^{p^{a-c}}$$
On the other hand, since $$z^{p^a}=y^{1-\sigma}$$
it comes $$z^{p^aA}=y^{(1-\sigma)A}$$
$$z^{p^aA}=\frac{N(y)}{y^{p^n}}$$
$$N(y)=z^{p^aA}y^{p^n}=r_0$$

\vspace{0,5cm}

Here, we assumed that $\mu_{p^a}\subset N(\ft)$, and then there exists $z\in\ft$ such that $N(z)=\xi$ where $\xi$ is a primitive root of unity of order $p^a$; since $N(z)^{p^a}=\xi^{p^a}=1$, the H90 theorem states that $H^1(G,\ft)=1$ and that we can find $y\in\ft$ such that $z^{p^a}=y^{1-\sigma}$. The lemma above then ensures that $N(y)=r_0$. We also consider $x\in\kt$ of order $p^n$ modulo $N(\ft)$.\vspace{0,5cm}

 {\bf lemma (3.2.2):} a convenient direct factor \vspace{0,5cm}
 
 When $\mu_K=\mu_F=\mu_{p^a}\in N(\ft)$, if $r_0$ denotes the radical of $K_1/K$ (unique modulo $K^{\times p}$),
 
 $$W_0=\langle\mu_K,r_0,x^{p^n}\rangle$$
 
 is a direct factor in $N(\ft)$.\vspace{0,5cm}
 
 Proof: \vspace{0,5cm}
 
 Since $\mu_K\subset W_0$, $W_0$ is a direct factor in $N(\ft)$ iff $W_0/\mu_K$ is a direct factor in $N(\ft)/\mu_K$, which has no torsion. Following (1.8), we assume that 
 
 $$\hspace{-5cm}(3.2.3)\hspace{3cm}N(t)^p=\xi^\alpha r_0^\beta x^{p^n\gamma}$$
 
  for some $\alpha,\,\beta, \gamma\in\zp$ and will prove that $p$ divides both $\beta$ and $\gamma$.
 
 At first, if $p\nmid\beta$ one has $r_0\in\mu_{p^a}K^{\times p}$, and then $F\supset K[\sqrt[p]{r_0}]=K[\mu_{p^{a+1}}]$, hence $\mu_F\neq\mu_K$, which contradicts our hypothesis. So $p\mid\beta$. It follows that $\xi^\alpha\in K^{\times p}$ and then $p\mid\alpha$; but then
  $$N(t)=\xi^{\alpha'} r_0^{\beta'} x^{p^{n-1}\gamma}$$
  
  for some $\alpha',\ \beta'\in\zp$ and then, since both $\xi$ and $r_0$ are in $N(\ft)$, $x^{p^{n-1}\gamma}\in N(\ft)$, so that $p\mid\gamma$.
 \vspace{0,5cm}

 This lemma allow us to consider the equality 
 $$N(\ft)=W_0\oplus S$$
 
 for some $S\simeq\zp^{d-1}$. We then can write $$S=\langle v_1,v_2,...,v_{d-1}\rangle$$
 with $$v_i=N(u_i),\ 1\le i\le d-1$$
 and then
 $$N(\ft)=\langle N(z),N(y),N(x), N(u_1),N(u_2),...,N(u_{d-1})\rangle$$
 so that, following the "N(L)=N(M)" lemma:
 $$\ft=\langle z,y,x,u_1,...,u_{d-1}\rangle$$

 Finally, let $$V=\frac{\langle X,Y,Z,U_1,...,U_{d-1}\rangle_{\zp[G]}}{\langle Z^{p^a}Y^{\sigma-1}\rangle_{\zp[G]}}=\frac{M}{D}$$
 
 be the formal space where, in $M$: $$\langle X,Y,Z,U_1,...,U_{d-1}\rangle=\langle X\rangle\oplus\langle Y\rangle\oplus\langle Z\rangle\oplus\langle U_1\rangle\oplus....\langle U_{d-1}\rangle\simeq \zp\oplus \zp[G]\oplus...\oplus\zp[G]$$ with $X^\sigma= X$ and, in $V$:
 
 $$\langle X\rangle\oplus\langle Y,Z\rangle\oplus\langle U_1\rangle\oplus....\langle U_{d-1}\rangle\simeq\zp\oplus\frac{\langle Y,Z\rangle}{\langle Z^{p^a}Y^{\sigma-1}\rangle}\oplus \zp[G]\oplus...\oplus\zp[G]$$
 which means that $X^\sigma=X$ and $Z^{p^a}=Y^{1-\sigma}$ generate all relations between $X,Y,Z,U_2...U_{d-1}$. We study in IV the structure of $W=\langle X\rangle\oplus\frac{\langle Y,Z\rangle}{\langle Z^{p^a}Y^{\sigma-1}\rangle}=\frac{\langle X,Y,Z\rangle}{\langle Z^{p^a}Y^{\sigma-1}\rangle}$, proving:
 $$H^0(G,W)=p^n$$
 $$H^1(G,W)=1$$
 $$\chi_W=\chi_{reg}+1$$
 and
 $$W_{tors}=\mu_{p^a}$$
 Then, $\ft$ is a quotient of $V$. But it has the same character
 $$\chi_{\ft}=\chi_V=d\chi_{reg}+1$$
 and has the same torsion, so that, noting $\phi$ the natural surjective map $\phi:V\rightarrow\ft,$ one has $\chi(\ker(\phi))=0$, then $ \ker(\phi)\subset V_{tors}=\ft_{tors}=\mu_{p^a}$, then
 $$\ker(\phi)={1}$$
 and $$\ft\simeq V$$
 
 \vspace{7cm}
 {\bf \large Case 3: $F=K[\mu_F]$: the cyclotomic case} \vspace{1cm}
 
 \begin{enumerate}
 	\item {\bf Discussion about $N(\mu_F)$:}\vx

 We first notice that $F=K[\mu_F]$ implies $a\ge 1$ because for $p\neq 2$, $[K[\mu_p]:K]$ is prime to $p$ and for $p=2$, $\mu_2\in K$. Also, $$\underset{\xi^p=1}{\Pi}\xi = \left\{ \begin{matrix}
 	1\ \text{ if }\ p\neq 2\\	-1\ \text{ if }\ p = 2
 \end{matrix}\right. $$

Hence, noting $\xi$ a primitive root of unity of order $p^{a+1}$,  $$N_{K[\mu_{p^{a+1}}]/K}(-\xi)=-\xi^p$$

When $p\neq 2$ or $a\ge 2$, the last equality propagates along $K[\mu_\infty]/K$ and one has both, noting $\xi_F$ a generator of $\mu_F=\mu_{p^{a+m}}$ and $\xi_K=\xi_F^{p^m}$ a generator of $\mu_K=\mu_{p^a}$:
$$N_{F/K}(-\xi_F)=-\xi_K$$
and
$$[K[\mu_{p^{a+m}}]:K]=p^m$$
However, when $p=2$ and $K\not\ni i$ this may fail, according to the following result:

\vspace{0,5cm}
{\bf lemma (3.3.1): } for any $\xi_F$ generating $\mu_F$ in the cyclic extension $F/K$, noting $[K[\mu_F]:K]=p^m$ and $K_m=K[\mu_F]$:

$$N_{{
	}K_m/K}(\xi_F)=\left\{\begin{matrix}
 & \xi_F^{p^m}&\text{if } p\neq 2\ \text{or}\ m=0\\
& & \\
 & -\xi_F^{2^m}&\text{if } p= 2,\  a\ge 2\text{ and }m>0\\
& & \\
&-1 &\text{if } m>0,\ p= 2, \ a=1\text{ and } K[\mu_\infty]/K \text{ procyclic}\\
& & \\
&1 &\text{if } \ p= 2, \ a=1\text{ and } K[\mu_\infty]/K \text{ not procyclic}
\end{matrix}
\right.$$

 \vspace{0,5cm}
 
 Proof: 
 \vspace{0,5cm}
 
 The case $p\neq 2$ or $a\ge 2$ as already been mentionned above, so we focus on the case $m>0,\ p=2$ and $a=1$. At first, $K[\mu_\infty]/K$ is procyclic iff $\gal (K[\mu_\infty]/K)$ has no torsion element. One has an isomorphism: 
 
  $$\phi:\left(\begin{matrix}
 	\gal(\mathbb{Q}_2[\mu_\infty]/\mathbb{Q}_2)&\longrightarrow&\mathbb{Z}_2^\star\\
 	\varphi:\xi\in\mu_\infty\mapsto\varphi(\xi)=\xi^k&\longmapsto& k
 \end{matrix}\right)$$
 
and $\gal(K[\mu_\infty]/K)$ is a subgroup of $\gal(\mathbb{Q}_2[\mu_\infty]/\mathbb{Q}_2)$. Then, $K[\mu_\infty]/K$ is procyclic iff 

$$\phi(\gal (K[\mu_\infty]/K))\ni -1$$

  \vx

\underline{Either $K[\mu_\infty]/K$ is procyclic:}
any open procyclic subgroup of $\gal(\mathbb{Q}_2[\mu_\infty]/\mathbb{Q}_2)$ is equal to 

$$I_l=\phi^{-1}(1+2^l\mathbb{Z}_2)$$

 or to
 
 $$J_l=\phi^{-1}(-1+2^l\mathbb{Z}_2)$$
 
 with $l\ge 2$ (for $l=1$, $I_1=J_1$ is not procyclic). $I_l$ fixes $i$ and then,
 
 $$\gal(K[\mu_\infty]/K)=J_l$$
 
 for some $l\ge 2$. Then, $\gal(K[\mu_\infty]/K(i))=J_l^2=\phi^{-1}((-1+2^l)^2)=\phi^{-1}(1-2^{l+1}+2^{2l})$, so that:
 
 $$\mu_{K(i)}=\mu_{2^{l+1}}$$
 
 and, since $\phi^{-1}(-1+2^l)$ generates $\gal(K(i)/K)$, it comes:
 
 $$N_{K(i)/K}(\xi_{2^{l+1}})=\xi_{2^{l+1}}^{1-1+2^l}=\xi_{2^{l+1}}^{2^l}=-1$$
 
 After that, we can apply in $K(m)/K(i)$ (replacing $a$ with $l+1$ and $m$ with $m-1$) the results found when $a\ge 2$ and then for any $k\in\mathbb{N}^\star$:

 $$[K[\mu_{2^{l+k}}]:K]=2^{k}\ \text{and}\ N_{K[\mu_{2^{l+k}}]/K}(\mu_{K[\mu_{2^{l+k}}]})=-1$$
 
 which proves that
 
 $$\hspace{-5cm}(3.3.2)\hspace{3cm}\mu_F=\mu_{2^{l+m}}\ \  \text{and}\ \ N_{K_m/K}(\xi_F)=-1$$
 \vx

 \underline{Either $K[\mu_\infty]/K$ is not procyclic:} 
 then, $\phi^{-1}(-1)\in\gal(K[\mu_\infty]/K)$ and generates $\gal(K(i)/K)$, so that:

$$\forall\xi\in\mu_{K(i)}, \ N_{K(i)/K}(\xi)=\xi\xi^{-1}=1$$ 

Moreover,
 
$$\gal(K[\mu_{2^{l+1}}]/K)\simeq \mathbb{Z}/2\mathbb{Z}\times\mathbb{Z}/2\mathbb{Z}$$

and then, since $K_m/K$ is cyclic:

$$K_m\subset K(i)$$

 Finally:

 $$\hspace{-5cm}(3.3.3)\hspace{3cm}m\le 1\ \text{and}\ N_{K_m/K}(\xi_F)=1$$

 \vspace{0,5cm}

 From now on, the case $m>0,\ p=2$, $a=1$, $K[\mu_\infty]/K$ procyclic will be referred as the "special case 1" and the case $m>0,\ p=2$, $a=1$, $K[\mu_\infty]/K$ not procyclic as the "special case 2" to avoid repetition of all these conditions.
 \vspace{0,5cm}
 
  \item {\bf the most general situation:}\vx
 
 We now assume that $F=K[\mu_F]$. It follows from the previous lemma that 
 $$N_{F/K}(\mu_F)=\mu_K$$
 except in the special case 2, studied in 3. At first, choosing as usual $x\in K$ of order $p^n$ in $\ft/N(\ft)$ and a generator $\xi_K$ of $\mu_K$, we prove that 
 $$W_0=\langle\mu_K, x^{p^n}\rangle$$
 is a direct factor in $N(\ft)$: to do so, we first note that it contains $\kt_{tors}$ and as such, is a direct factor in $N(\ft)$ iff $\langle x^{p^n}\rangle$ is one of $N(\ft)/\mu_K$; assuming that $$N(t^p)=\xi_K^\alpha x^{p^n\beta}$$
 for some $\alpha,\ \beta\in\zp$, we first deduce that $\xi_K^\alpha\in K^{\times p}$, then $p\mid\alpha=p\alpha '$ and then,
 $$N(t/\xi^\alpha)^p=x^{p^n\beta}$$
 so that 
 $$N(t/\xi^\alpha)=\xi x^{p^{n-1}\beta}$$
 for some $\xi\in\mu_p$, which is a norm in $F/K$. Hence, $x^{p^{n-1}}$ is a norm as well, and $p\mid\beta$, as expected.
 
 It's now legit to choose $v_1,...,v_d\in N(\ft)$ such that
 
 $$N(\ft)=\langle \xi_K,x^{p^n},v_1,...,v_d\rangle_{\zp}$$
 Applying the "$N(M)=N(L)$" lemma and choosing $u_1,...,u_d\in\ft$ such that $N(u_i)=v_i$ for $1\le i\le d$ implies 
 
 $$\ft=\langle\xi_F,x,u_1,...,u_d\rangle_{\zp[G]}$$
 and then finally, after comparing characters and using the same arguments as in Case 2, one has
 
 $$\ft\simeq \mu_F\times\zp\times \zp[G]^d$$
 
 with $G$ acting $\omega$-isotypically on the first factor, trivially on the second, the last being a free module generated by $d$ elements.
 
 \vx
 
 \item {\bf  special case}\vx

Special case 1 was studied along the general situation in 2. What happens in special case 2 ?

(3.3.3) proved that $m\le 1$. Nothing really interesting has to be said if $m=0$, so we assume $m=1$ and then $F=K(i)$. We first note $l\ge 2$ the integer such that:

$$|\mu_F|=2^l$$

and $\xi_F$ denotes a generator of $\mu_F$. At last, we define

$$t=1-\xi_F$$

and make this computation:

$$t^{1-\sigma}=\frac{1-\xi_F}{1-\xi_F^{-1}}=\frac{1-\xi_F}{(\xi_F-1)\xi_F^{-1}}=-\xi_F$$

and then,

$$\hspace{-5cm}(3.3.4)\hspace{3cm}t^{2^{l-1}(1-\sigma)}=-1$$

our fundamental relation will derive directly from this, but at first we have to wonder if $-1$ is a norm in $F/K$ or not:

\vx

{\bf lemma (3.3.5):}

\vx

In special case 2, -1 is a norm in $K(i)/K$ iff $\dim_{\mathbb{Q}_2}(K)$ is even.

\vx

Proof: This is a direct consequence on the compatibility between the norm and restriction applications, see (1.2) replacing in the diagram $K$ with $\mathbb{Q}_2$ and $L$ with $K$.\vx

\begin{enumerate}
	\item \underline{$\dim_{\mathbb{Q}_2}(K)$ is odd:}
		
		\vx
		
		Here, $-1$ is not a norm in $F/K$ with $F=K(i)$. One has, noting $\theta=\frac{2\pi}{2^l}$ (see also (3.6.4)):
		
		$$N(t)=2-2\cos(\theta)$$

		which is not a square in $K$ (proved in lemma (3.6.9)). Moreover, $N(\ft)\not\ni -1$ and then, $N(\ft)$ has no torsion. Then,
		
		$$W_0=\langle t\rangle_{\zp}$$
		
		is a direct factor in $N(\ft)\simeq\mathbb{Z}_2^{d+1}$. We then choose $S$ such that 
		
		$$N(\ft)=W_0\oplus S$$
		
		and a family $v_1,...,v_d$ generating $S$. We then write $v_i=N(u_i)$ for some $u_i\in\ft$ and $1\le i\le d$. The "N(L)=N(M)" lemma proves that:
		
		$$\ft=\langle t, u_1,...,u_d\rangle_{\mathbb{Z}_2[G]}$$

		Finally, let $$V=\frac{\langle T,U_1,...,U_{d}\rangle_{\zd[G]}}{\langle T^{2^l(1-\sigma)}\rangle_{\zd[G]}}=\frac{M}{D}$$
		
		be the formal space where $M$ is a free $\zd[G]$-module generated by the $d+1$ elements $T,U_1,...,U_d$ and, in $V$:
		
		$$\langle T\rangle\oplus\langle U_1\rangle\oplus....\langle U_{d}\rangle\simeq\frac{\langle T\rangle}{\langle T^{2^l(1-\sigma)}\rangle}\oplus \zd[G]\oplus...\oplus\zd[G]$$
	 
	 \vx
	 
	  We now study $W_{tors}$, where

	   $$W=\frac{\langle T\rangle}{\langle T^{2^l(1-\sigma)}\rangle}$$\vx

 If $T^\alpha\in W_{tors}$ has order $p^k$ in $W$, then in $\zd[G]$:

$$2^k\alpha =2^l(1-\sigma)\lambda$$

for some $\lambda\in\zd[G]$. We then make the euclidean division of $\lambda$ by $N$, seeing $\lambda$ as a polynomial of degree $i<2^n=2$:

$$\lambda=N\times Q+R$$

with $R\in\zd[X]$ of degree $d<d^\circ N=1$, hence $R$ is a constant and

	$$2^k\alpha =2^l(1-\sigma)(NQ+R)=2^l(1-\sigma)R$$
	
	Then, $2^k\mid 2^lR=2^kC$ for some $C\in\zd$ and $\alpha=C(1-\sigma)$. This proves that $T^{1-\sigma}$ generates $W_{tors}$, which order is $2^l$.\vx

	\vx

		To conclude, we first note that $\ft$ is a quotient of $V$. But it has the same character
		$$\chi_{\ft}=\chi_V=d\chi_{reg}+1$$
		(because $\chi(W)=1$) and has the same number of torsion elements, so that, noting $\phi$ the natural surjective map $\phi:V\rightarrow\ft,$ one has $\chi(\ker(\phi))=0$, then $ \ker(\phi)\subset V_{tors}$. Hence, $\ker(\phi)\neq{1}$ would implies that $|\ft_{tors}|<|V_{tors}|$ which is not possible. 
		Finally:
		
		$$\ft\simeq V$$ 
		
		{\bf Remark:} no $x$ generating $\kt/N(\ft)$ appears here. However, $x=-1$ would be convenient, and is contained in $\langle t \rangle =W$. \vx

  	\item \underline{$\dim_{\mathbb{Q}_2}(K)$ is even:}
  
  \vx
  
  Here, $-1$ is a norm in $F/K$ and $F=K(i)$. We then choose $y$ such that 
  
  $$N(y)=y^{1+\sigma}=-1$$
  
  and it comes from (3.3.4):
  
  $$\hspace{-5cm}(3.3.6)\hspace{3cm}t^{2^{l-1}(1-\sigma)}=y^{1+\sigma}$$
  
  (3.3.6) is our fundamental relation. At last, we choose $x\in\ft$ generating $\kt/N(\ft)$.
  
  \vx
  {\bf lemma (3.3.7):} $W_0=\langle x^2,N(y),N(t)\rangle$ is a direct factor in $N(\ft)$.
  \vx
  
  Proof: $W_0\supset\kt_{tors}=N(\ft)_{tors}=\mu_2$ so that $W_0$ is a direct factor of $N(\ft)$ iff $W_0/\mu_2$ is one of $N(\ft)/\mu_2$. Assuming
  
  $$N(f^2)=(-1)^k x^{2\alpha} N(t)^\beta N(y)^\gamma$$
  
  we then must prove that $N(f)\in \mu_2W_0$. We begin rewriting the previous equality, using $N(y)=-1$:
  
  $$\hspace{-5cm}(3.3.7)\hspace{3cm}N(f^2)=\pm\ x^{2\alpha} N(t)^\beta $$
  
  It follows that $\pm \ N(t)^\beta$ is a square in $\kt$. lemma (3.6.9) then proves that $2|\beta=2\beta '$ and then, (3.3.7) becomes:
  
  $$N\left(\frac{f}{t^{\beta '}}\right)^2=\pm\ x^{2\alpha} $$
  
  and then,
  
  $$N\left(\frac{f}{t^{\beta '}}\right)=\pm\ x^{\alpha} $$
  
  Whatever the sign is, $-1$ is a norm so that $x^\alpha$ as well, and then $2\mid\alpha=2\alpha '$, and finally:
  
  $$N(f)= \pm\ x^{2\alpha '} N(t)^{\beta '} \in \mu_2 W_0$$
  
  This proved as well (see 1.8) that $\dim_{\mathbb{F}_2}\langle x^2,N(t)\rangle=2$ in $\frac{N(\ft)}{\mu_2N(F^{\times })^2}$, and then $\chi(W_0)=2\times 1$. We go on choosing $S\simeq \zd^{d-1}$ such that 
  
  $$N(\ft)=W_0\oplus S$$
  
  and a family $v_1,...,v_{d-1}$ generating $S$. We then write $v_i=N(u_i)$ for some $u_i\in\ft$ and $1\le i\le d-1$. The "N(L)=N(M)" lemma proves that:
  
  $$\ft=\langle x,y,t, u_1,...,u_{d-1}\rangle_{\mathbb{Z}_2[G]}$$

  Finally, let $$V=\frac{\langle X,Y,T,U_1,...,U_{d-1}\rangle_{\zd[G]}}{\langle T^{2^{l-1}(1-\sigma)}Y^{-1-\sigma}\rangle_{\zd[G]}}=\frac{M}{D}$$
  
  be the formal space where the equality $X^\sigma=X$ both holds in $M$ and $V$, and where in $V$:
  
  $$\langle X,Y,T\rangle\oplus\langle U_1\rangle\oplus....\langle U_{d-1}\rangle\simeq\zd\oplus\frac{\langle Y,T\rangle}{\langle T^{2^{l-1}(1-\sigma)}Y^{-1-\sigma}\rangle}\oplus \zd[G]\oplus...\oplus\zd[G]$$
  
  \vx
  
  We now have a look at $$W=\frac{\langle Y,T\rangle}{\langle T^{2^{l-1}(1-\sigma)}Y^{-1-\sigma}\rangle}$$\vx
  
 The main point is to determine $W_{tors}$. To do so, we assume that $w=T^\alpha Y^\beta\in W_{tors}$ has order $p^k$ in $W$, for some $k\ge 1$; then, in $\zd[G]$:
 
 $$\left\{\begin{matrix}
 	2^k\alpha=2^{l-1}(1-\sigma)\lambda\\
 	2^k\beta=-(1+\sigma)\lambda
 \end{matrix}\right.$$
 
 for some $\lambda\in\zd[G]$. But then, $(1-\sigma)\mid\alpha$ and $(1+\sigma)\mid\beta$. It follows that 
 
 $$w\in\langle T^{1-\sigma}, N(Y)\rangle_{\zd[G]}$$
 
 and in $W$, $N(Y)=T^{2^{l-1}(1-\sigma)}$ so that 
 
 $$w\in\langle T^{1-\sigma}\rangle_{\zd[G]}$$
 
 Moreover, $T^{2^{k}(1-\sigma)}=1$ in $W$ iff 
 
 $$\left\{\begin{matrix}
 	2^k(1-\sigma)=2^{l-1}(1-\sigma)\lambda\\
 	0=-(1+\sigma)\lambda
 \end{matrix}\right.$$

We can see $\lambda$ as a polynomial of degree $d<p^n=2$ then we can write $\lambda=r+s\sigma$ and the second equation holds iff $s=-r$ and $\lambda=r-r\sigma$. The first then becomes, since $(1-\sigma)^2=2-2\sigma\in\zd[G]$:

$$2^k(1-\sigma)=2^{l-1}(1-\sigma)(r-r\sigma)=2^{l}r(1-\sigma)$$
 
 This implies $k\ge l$; conversely in $W$:
 
 $$T^{2^{l}(1-\sigma)}=N(Y)^2=N(N(Y))=N(T^{2^{l-1}(1-\sigma)})=1$$
 
 This proves that:
 
 $$|W_{tors}|=2^l$$
 
 as expected.\vx
 
 The next step is computation of $\chi(W)$. In $\mathbb{Q}_2\underset{\zd}{\otimes}W$, 
 
 $$T^{2^{l-1}(1-\sigma)}=N(Y)$$
 
 implies $N(N(Y))=1$ and then $N(Y)=1$, and as well $T^{1-\sigma}=1$. Then,
 
 $$\chi(W)=\chi(\langle T\rangle)+\chi(\langle Y\rangle)=1+(\chi_{reg}-1)=\chi_{reg}$$
 
 and then, 
 
 $$\chi(V)=d\chi_{reg}+1=\chi(\ft)$$
 
 as expected. Then, $V$ and $\ft$ share the same character and the same number of torsion elements. $\ft$ is a quotient of $V$ and using the usual arguments:
 
 $$\ft\simeq V$$

\end{enumerate}

  \end{enumerate}

 \vspace{1cm}
 {\bf \large Case 4: $\mu_F=\mu_K$} \vspace{1cm}
 
 This situation is similar to Case 2, but more general since we no more suppose that $\mu_K\in N(\ft)$. We'll then use the parameter $b$ as defined at the beginning of III. 
 
 \vspace{0,5cm}
 
 As usual we first choose $x$ generating $\ft/N(\ft)$. We then notice that the order of $\xi_K$ modulo $N(\ft)$ being $p^b$, one has
 
 $$\xi_K\equiv x^{\alpha p^{n-b}}\ \mod N(\ft)$$
 
 for some $\alpha\in\zp^\times$. Replacing $x$ with $x^
{-\alpha}$ generating as well $\kt/N(\ft)$, it comes:

$$\hspace{-5cm}(3.4.1) \hspace{3cm} \xi_K x^{p^{n-b}}=N(s)$$
for some $s\in\ft$. It follows that 

$$N\left(s^{p^a}/x^{p^{a-b}}\right)=(\xi_K x^{p^{n-b}})^{p^a}x^{p^{n+a-b}}=1$$

which allows to write:

$$\hspace{-5cm}(3.4.2) \hspace{3cm}\frac{s^{p^a}}{x^{p^{a-b}}}=t^{1-\sigma}$$

for some $t\in\ft$. (3.4.2) is our fundamental relation. We then define 

$$r=\left\{\begin{matrix}
	s^A t^{p^{n-a}} \text{ if }\ n\ge a\\
	s^{p^{a-n}A} t \text{ if }\ n\le a
\end{matrix}
\right.$$

Denoting $c=\min (n,a)$, it comes:

$$r=s^{p^{a-c}A} t^{p^{n-c}} $$
and 
$$r^{1-\sigma}=s^{p^{a-c}A(1-\sigma)} t^{p^{n-c}(1-\sigma)} $$
Using (1.7) along with (3.4.2):
$$r^{1-\sigma}=s^{p^{a-c}(N-p^n)} \left(\frac{s^{p^a}}{x^{p^{a-b}}}\right)^{p^{n-c}} $$

$$r^{1-\sigma}=\frac{N(s)^{p^{a-c}}}{s^{p^{n+a-c}}}\left(\frac{s^{p^a}}{x^{p^{a-b}}}\right)^{p^{n-c}}=\xi_K^{p^{a-c}} $$

which is a root of unity of order $p^c$. Then, 
$$r_0=r^{p^c}\in \kt$$

is a generator of $Rad(K_c/K)$ - note that $c>0$. On the other hand, 

$$t^{(1-\sigma)A}=\frac{s^{p^a A}}{x^{p^{a-b}A}}=\frac{N(t)}{t^{p^n}}$$
so that

$$N(t)=\frac{s^{p^a A}t^{p^n}}{x^{p^{a-b}A(1)}}$$

where $A(1)=\frac{p^n (p^n-1)}{2}$ is the sum of coefficients of $A$.

$$\hspace{-5cm}(3.4.3)\hspace{3cm}N(t)=\frac{r_0}{x^{p^{n+a-b}(p^n-1)/2}}$$

We then use the following result:

\vx

 {\bf lemma (3.4.4)}: $W_0=\langle N(s),N(x),N(t)\rangle $ is a direct factor in $N(\ft)$ for a convenient $x$.
 
 \vx
 
 Proof: We first note that 
 $$\xi_K^{p^b}=\frac{N(s^{p^b})}{N(x)}\in W_0$$
 so that $N(\ft)_{tors}\in W_0$, hence $W_0$ is a direct factor in $N(\ft)$ iff $\frac{W_0}{\mu_K\cap N(\ft)}$ is one of $\frac{N(\ft)}{\mu_K\cap N(\ft)}$.
 Using (1.8), we assume that
 
 $$N(f^p)=N(s)^\alpha N(x)^\beta N(t)^\gamma$$
 
 and have to prove that $N(f)\in W_0$. It comes:
 
 $$N(f^p)=\xi_K^\alpha x^{\alpha p^{n-b}+\beta p^n-\gamma p^{a-b+n}(p^n-1)/2}\ r_0^\gamma$$
 
 At first, notice that $p^{n-b}\mid (\alpha p^{n-b}+\beta p^n-\gamma p^{a-b+n}(p^n-1)/2)$, even if $p=2$ (because $a>0$).

 If $n>b$, it follows that $r_0^\gamma=\xi_K^{-\alpha}\delta^p$ for some $\delta\in \kt$, then $p\mid \gamma$ or $\mu_{p^{a+1}}\subset F$; since $\mu_K=\mu_F$, the last option is not possible and $p\mid \gamma$.
 
 If $n=b$, the first thing to say is that choosing $x$ at random was definitely not the best choice, which was $x=\xi_K$; however, to avoid the study of more and more cases, we'll prove that lemma (3.4.4) remains valid. At first, we remember that when $x$ generates $\kt/N(\ft)$, $xN(y)$ as well, for any $y\in\ft$: the choice of $x$ is modulo $N(\ft)$. Then, we can choose $x$ such that 
 
 $$\hspace{-5cm}(3.4.5)\hspace{3cm}\dim_{\mathbb{F}_p}\langle\xi_K,x,r_0\rangle = 3$$ in $\kt/K^{\times p}$ because if not, we have both $x\in\langle\xi_K,r_0\rangle$ and $N(\ft)\subset \langle\xi_K,r_0\rangle$ modulo $K^{\times p}$, and since $dim_{\mathbb{F}_p}N(\ft)=d+1$ this would imply $d=1$, $N(\ft)= \langle\xi_K,r_0\rangle\ni x\mod N(\ft)$ and then $x\in N(\ft)K^{\times p}$ cannot generate $\kt/N(\ft)$, which is a contradiction. 
 
 So, we proved we can choose $x$ satisfying (3.4.5) and from now on, will suppose it. But it then becomes obvious that in the penultimate equality, $p$ divides both $\alpha=p\alpha'$ and $\gamma=p\gamma'$.
 
  Then, 
  $$N\left(\frac{f}{t^{\gamma '}s^{\alpha '}}\right)^p= x^{\beta p^n}$$
  so that 
  $$N\left(\frac{f}{t^{\gamma '}s^{\alpha '}}\right)=\xi x^{\beta p^{n-1}}$$
  
  for some $\xi\in\mu_p$. Noting $\xi=\xi_K^\lambda$, it comes
  
  $$N\left(\frac{f}{t^{\gamma '}s^{\alpha '}}\right)=\xi_K^\lambda x^{\beta p^{n-1}}=N(s^\lambda)x^{\beta p^{n-1}-\lambda p^{n-b}}$$
  
   so finally
   
    $$N\left(\frac{f}{t^{\gamma '}s^{\alpha '+\lambda }}\right)=\xi_K^\lambda x^{\beta p^{n-1}}=x^{\beta p^{n-1}-\lambda p^{n-b}}$$
    
    The right term of this equality is a norm and can then be written $N(x^\mu)$ for some $\mu\in\zp$, hence $$N(f)=N(x^\mu s^{\alpha ' +\lambda} t^{\gamma '})\in W_0$$
    
    and the lemma is proved.
    
    \vx

    It follows from (3.4.5) that $\dim_{\mathbb{Q}_p}(\qp\underset{\zp}{\otimes}W_0)=2$ and then, writing 
    
    $$N(\ft)=W_0\oplus S$$
    
    one has, since $N(\ft)_{tors}\subset W_0$: 
    
    $$S\simeq\zp^{d-1}$$
    
    We then can write $$S=\langle v_1,v_2,...,v_{d-1}\rangle$$
    with $$v_i=N(u_i),\ 1\le i\le d-1$$
    and then
    $$N(\ft)=\langle N(s),N(t),N(x), N(u_1),...,N(u_{d-1})\rangle$$
    so that, following the "N(L)=N(M)" lemma:
    $$\ft=\langle s,t,x,u_1,...,u_{d-1}\rangle$$

    Finally, let $$V=\frac{\langle X,S,T,U_1,...,U_{d-1}\rangle_{\z2[G]}}{\langle S^{p^a}X^{-p^{a-b}}T^{\sigma-1}\rangle_{\z2[G]}}=\frac{M}{D}$$
    
    be the formal space where, in $M$: $$\langle X,S,T,U_1,...,U_{d-1}\rangle=\langle X\rangle\oplus\langle S\rangle\oplus \langle T\rangle\oplus\langle U_1\rangle\oplus....\langle U_{d-1}\rangle  \simeq  \zp\oplus \zp[G]\oplus...\oplus\zp[G]$$
    
    and in $V$:
    
    $$\langle X,S,T\rangle\oplus\langle U_1\rangle\oplus....\langle U_{d-1}\rangle   \simeq    \frac{\langle X,S,T\rangle}{\langle S^{ p^a}X^{-p^{a-b}}T^{\sigma-1}\rangle}\oplus \zp[G]\oplus...\oplus\zp[G]$$
    which means that $X^\sigma=X$ and $\frac{S^{p^a}}{x^{- p^{a-b}}}=T^{1-\sigma}$ generate all relations between $X,S,T,U_1...U_{d-1}$. We study in IV the structure of $W=\frac{\langle X,S,T\rangle}{\langle S^{p^a}X^{-p^{a-b}}T^{\sigma-1}\rangle}$, proving:
    $$H^0(G,W)=p^n$$
    $$H^1(G,W)=1$$
    $$\chi_W=\chi_{reg}+1$$
    and
    $$W_{tors}=\mu_F=\mu_K$$
    Then, $\ft$ is a quotient of $V$. But it has the same character
    $$\chi_{\ft}=\chi_V=d\chi_{reg}+1$$
    and has the same torsion, so that, using the same arguments as at the end of Case 2:
     $$\ft\simeq V$$

\vspace{1cm}
{\bf \large Case 5: the mixed case} \vspace{1cm}

We now study the most general case, which is a mix between the cyclotomic case and Case 4. The situation here is not absolutely general: we assume $\mu_F\neq 1$ (otherwise Case 1 applies). Another difficulty is that we don't suppose that $\mu_K\subset N(\ft)$. We then have to use the parameters $b$ and $m$ defined at the beginning of III, and will assume that $m\neq 0$ (otherwise Case 4 applies) and that $n>m$ (otherwise the cyclotomic case applies). Since $\mu_F\neq 1$, $\mu_K\neq 1$ (because if $p\neq 2, [K[\mu_p]:K]$ is prime to $p$ and if $p=2$, $K\supset\mu_2$ anyway), so $a\ge 1$. Moreover, the special case 2 is excluded here (Case 6 will apply), but not the special case 1 (in which the parameter $l\ge 2$ such that $\mu_{K(i)}=\mu_{2^l}$ will be used). As above, we choose $x\in\ft$ of order $p^n$ modulo $N(\ft)$.

  As special case 2 is excluded, lemma (3.3.1) states that
  
  $$N_{K_m/K}(\mu_F)=\mu_K$$
  
  and then, after choosing a generator $\xi_F$ of $\mu_F$, we define
   
   $$\hspace{-5cm} (3.5.1)\hspace{3cm}\xi_K=N_{K_m/K}(\xi_F)$$

At this point, we note that the order of $x$ in $\frac{K_m^\times}{N_{F/K_m}(\ft)}$ is $p^{n-m}$ (it can't be more because $p^{n-m}=[F/K_m]$, and $x^{p^{n-m-1}}=N_{F/K_m}(t)$ would imply $N_{F/K}(t)=x^{p^{n-1}}$ and would be a contradiction). Then, $x$ generates as well $\frac{K_m^\times}{N_{F/K_m}(\ft)}$ and as such it's legitimate to write, $N_m$ denoting $N_{F/K_m}$:

$$N_m(s)=\xi_Fx^\alpha$$
for some $s\in\ft$ and some $\alpha\in\zp$. Taking $N_{K_m/K}$ of this equality leads to 

$$N_{F/K}(s)=\xi_Kx^{\alpha p^m}$$
and since the order of $\xi_K$ modulo $N(\ft)$ is $p^b$ and $\xi_K\equiv x^{-\alpha p^m}\ \mod N(\ft)$  it follows

$$\hspace{-5cm} (3.5.2)\hspace{3cm}n-b\ge m$$
and $\alpha=\beta p^{n-b-m}$ for some $\beta\in\zp^\times$. We then replace $x$ with $x^\beta$, convenient as well, so that 2 of the last 3 eqalities become:

$$\hspace{-5cm} (3.5.3)\hspace{3cm}N_m(s)=\xi_Fx^{p^{n-m-b}}$$

and then,

$$\hspace{-5cm} (3.5.4)\hspace{3cm}N_{F/K}(s)=\xi_Kx^{p^{n-b}}$$

We now have to introduce $k_\sigma\in\zp$ such that:

$$\hspace{-5cm} (3.5.5)\hspace{3cm}\forall\xi\in\mu_F:\xi^\sigma=\xi^{k_\sigma}$$

One has:

 $$k_\sigma=1+\kappa p^a$$
 
 for some $\kappa\in\zp^\times$ (because $\sigma$ has to fix any element of $\mu_{p^a}$ ant not all of $\mu_{p^{a+1}}$). Except in the special cases, you may assume that $\kappa=1$, eventually replacing $\sigma$ with $\sigma^\alpha$ for some $\alpha\in\zp^\times$. However we won't suppose $\kappa=1$.

Since $$N_m(s)=\xi_Fx^{p^{n-m-b}}$$
one has:
$$N_m(s)^{-\sigma+k_\sigma}=x^{(k_\sigma -1)p^{n-m-b}}=x^{\kappa p^{a+n-m-b}}$$
 and then:

$$ N_m\left(\frac{s^{-\sigma+1+\kappa p^a}}{x^{\kappa p^{a-b}}}\right)=1$$

applying H90 in $F/K_m$ of which the Galois group is generated by $\sigma^{p^m}$:

$$\hspace{-5cm}(3.5.6)\hspace{3cm}\frac{s^{-\sigma+1+\kappa p^a}}{x^{\kappa p^{a-b}}}=t^{1-\sigma^{p^m}}$$ 

for some $t\in \ft$.  
 
(3.5.6) is our fundamental relation, playing the very same role as $z^{p^a}=y^{1-\sigma}$ in Case 2.\vspace{0,5cm}

The next step is computation of $N_{F/K}(t)$. To do so, we'll use the 2 following lemma: 
\vspace{0,5cm}

{\bf lemma (3.5.7):} computing $N(y)$ when $z^{p^a}=x_0y^{1-\sigma}, N(z)^{p^{a-c}}=\xi x_0^{p^{n-c}},\ x_0\in K$, with $\xi$ of order $p^c$\ (with $c=\min (n,a)$)\vspace{0,5cm}

Let $F/K$ be a finite cyclic extension of fields (local or not) of degree $p^n$, with $\mu_{p^a}\subset K$ and $c=\min (n,a)$. Then, the relation $$z^{p^a}=x_0y^{1-\sigma}$$ for $z$ satisfying $N(z)^{p^{a-c}}=\xi x_0^{p^{n-c}}$ with $\xi$ of order $p^c$, implies that, one has:

$$r=z^{p^{a-c}A}y^{p^{n-c}}$$

is such that $$r^{(1-\sigma)}=\xi$$
so that $$r_0=r^{p^c}\in K$$
is the radical of the unique extension of degree $p^c$ contained in $F$. Moreover, one has:
$$N(y)=\frac{r_0}{x_0^{p^n(p^n-1)/2}}$$

\vspace{0,5cm}

Proof:\vspace{0,5cm}

Like in lemma (3.2.1) with a few differences:\vspace{0,5cm}

$$r=z^{p^{a-c}A}y^{p^{n-c}}$$
$$r^{(1-\sigma)}=z^{p^{a-c}A(1-\sigma)}y^{p^{n-c}(1-\sigma)}$$
$$r^{(1-\sigma)}=\frac{N(z)^{p^{a-c}}}{z^{p^{n+a-c}}}\frac{z^{p^{a+n-c}}}{x_0^{p^{n-c}}}$$
$$r^{(1-\sigma)}=\frac{N(z)^{p^{a-c}}}{x_0^{p^{n-c}}}=\frac{\xi x_0^{p^{n-c}}}{x_0^{p^{n-c}}}=\xi$$
of order $p^c$. On the other hand, since $$z^{p^a}=x_0y^{1-\sigma}$$
it comes $$z^{p^aA}=x_0^Ay^{(1-\sigma)A}$$
$$z^{p^aA}=\frac{x_0^AN(y)}{y^{p^n}}$$
$$N(y)=\frac{z^{p^aA}y^{p^n}}{x_0^A}=\frac{r_0}{x_0^A}=\frac{r_0}{x_0^{p^n(p^n-1)/2}}$$
because $A(1)=\frac{p^n(p^n-1)}{2}$ is the sum of coefficients of $A$ and $x_0^\sigma=x_0$ (note that for any $P\in\zp[G]=\zp[X]/(X^{p^n}-1)$, and for any $\xi\in\mu_{p^n}$ we can define $P(\xi)$, which legitimates the usage of the notation $A(1)$).

\vspace{0,5cm}

 We also need the following result before looking for a convenient direct factor in $N(\ft)$:
 \vspace{0,5cm}
 
 {\bf lemma (3.5.8):} Noting $c=\min(n-m,a)$, one has: $Rad(K_{m+c}/K_m)=\xi_F\delta_0$ for some $\delta_0\in\kt$.
 
 \vspace{0,5cm}
 
 Proof: At first, note that $Rad(K_{m+c}/K_m)$ can't be chosen in $\kt$ because the compositum of $F$ and $K[\sqrt[p]{\delta_0}]$ is either $F$ or not cyclic; also, $\mu_F=\mu_{K_m}$ by definition, and $G$ acts $\omega$-isotypically on $\langle r_m\rangle$, where $\omega$ is the character of the action of $G$ on $\mu_F$ (associated with $k_\sigma$). An $\omega$-isotypical action on elements of order $o\le p^a$ is simply the trivial action, and then 
 
 $$r_m\in\left(\frac{K_m^{\times  }}{K_m^{\times p^c  }}\right)^G$$
 
 so that $r_m^{1-\sigma }=y^{p^c}$
 for some $y\in K_m^{\times}$,
 hence $N_{Km/K}(y)^{p^c}=1$. If $N_{Km/K}(y)^{p^{c-1}}=1$, we can write $y^{p^{c-1}}=f^{1-\sigma}$ and then $r_m^{1-\sigma}=f^{(1-\sigma)p}$. It would follow $r_m=f^p\gamma_0$ for some $\gamma_0\in\kt$ and then $K_{m+1}=K_m[\sqrt[p]{\gamma_0}]$ which is impossible, as discussed above. So, $N_{K_m/K}(y)$ is a root of unity of order $p^c$ like $\xi_K^{p^{a-c}}$ and then, eventually replacing $y$ with $y^\beta$ with $\beta\in\zp^\times$:
 
 $$N_{Km/K}(y)=\xi_K^{p^{a-c}}=N_{Km/K}(\xi_F^{p^{a-c}})$$
 
 so that $y=\xi_F^{p^{a-c}} w^{1-\sigma}$ for some $w\in K_m^\times$ and then, 
 
 $$r_m^{1-\sigma }=y^{p^c}=\xi_F^{p^a} w^{(1-\sigma)p^c}=\xi_F^{-\kappa^{-1}(1-\sigma)}w^{(1-\sigma)p^c}$$
 since $\xi_F^\sigma=\xi_F^{1+\kappa p^a}$. Finally, as expected:
 
 $$r_m=\xi_F^{-\kappa^{-1}}w^{p^c}\delta_0\equiv \xi_F^{-\kappa^{-1}}\delta_0 \ \ \ mod\ K_m^{\times p^c}$$
 for some $\delta_0\in\kt$. Of course, we can replace it with $r_m^{-\kappa}$ and rename $\delta_0^{-\kappa}$ as $\delta_0$, so that the lemma (3.5.8) is proved.

 \vspace{0,5cm}
 
 {\bf remark (3.5.9):} This proves that finally, altough $K_{m+c}/K$ is not Kummer, you can write:
 $$K_{m+c}=K_m\left[\sqrt[p^c]{\xi_F\delta_0}\right]=K\left[\sqrt[p^{c+m}]{\xi_K\delta_0^{p^m}}\right]$$
 The last equality fails in the special case 1, where $c=a=1$ and $K\left[\sqrt[2^{m+l-1}]{-\delta_0^{2^{ ^{m+l-2}}}}\right]$ is the correct formula, but it's not of high interest.

 \vx

 Now, we're looking for such a "radical" 
 
 $$R=r^{p^{m+c}}\equiv \left(\xi_K\delta_0^{p^m}\right)^\tau\ \mod K^{\times p^{m+c}}$$ 
 
 as described in the previous remark with $\tau\in\zp^\times$ and $r$ satisfying $r^{1-\sigma}=\xi$ a root of unity of order $p^{m+c}$ (or $2^{m+l-1}$ in the special case 1). Note that when $\xi=r^{1-\sigma}$ is fixed, $R$ is unique in $\kt/K^{\times p^{m+c}}$ (because if $r^{1-\sigma}=r'^{1-\sigma}=\xi$, then $r=r'\gamma_0$ for some $\gamma_0\in K$ and then, $R=R'\gamma_0^{p^{m+c}}$)  and that $c=\min(n-m,a)>0$; however here we won't have to pay attention to $\xi$, then $R$ is known up to some invertible power we denote $\tau$. 
 
 A generic method consists in applying lemma (3.5.8) in both $K_c/K$ and $K_{m+c}/K_m$ and then compare results, but we offer the following shortcut, using 
 
 $$S_m=1+\sigma+\sigma^2+...+\sigma^{p^m}\ \text{ and }\ A_m=\sum_{k=0}^{p^{n-m}-1}k\sigma^{kp^m}$$ satisfying $$1-\sigma^{p^m}=(1-\sigma)S_m\ \text{ and }\ A_m(1-\sigma^{p^m})=N_m-p^{n-m}\ \ :$$
 
 Let $r\in K_{m+c}$ be the number defined by
 
 $$\hspace{-5cm}(3.5.10)\hspace{3cm}r=s^{\kappa S_mA_mp^{a-c}} \left(\frac{t^{S_m}}{s}\right)^{p^{n-m-c}}$$
 
 Then, one has:
 
 $$r^{1-\sigma}=s^{\kappa (1-\sigma^{p^m})A_mp^{a-c}} \left(\frac{t^{S_m}}{s}\right)^{(1-\sigma)p^{n-m-c}}$$

$$r^{1-\sigma}=s^{\kappa (N_m-p^{n-m})p^{a-c}} \left(\frac{t^{S_m}}{s}\right)^{(1-\sigma)p^{n-m-c}}$$

On the other hand, one has from (3.5.6):

$$\frac{s^{\kappa p^a}}{x^{\kappa p^{a-b}}}=\left(\frac{t^{S_m}}{s}\right)^{1-\sigma}$$

so that the penultimate equality becomes:

$$r^{1-\sigma}=s^{\kappa (N_m-p^{n-m})p^{a-c}} \left(\frac{s^{\kappa p^a}}{x^{\kappa p^{a-b}}}\right)^{p^{n-m-c}}$$

$$r^{1-\sigma}=\frac{N_m(s^{\kappa p^{a-c}})}{s^{\kappa p^{n-m+a-c}}} \left(\frac{s^{\kappa p^a}}{x^{\kappa p^{a-b}}}\right)^{p^{n-m-c}}$$

$$r^{1-\sigma}=\frac{N_m(s)^{\kappa p^{a-c}}}{x^{\kappa p^{n-m-c+a-b}}} $$

then using (3.5.3):

$$r^{1-\sigma}=\frac{(\xi_F x^{p^{n-m-b}})^{\kappa p^{a-c}}}{x^{\kappa p^{n-m-c+a-b}}} $$

$$\hspace{-5cm}(3.5.11)\hspace{3cm}r^{1-\sigma}=\xi_F^{\kappa p^{a-c}}$$

is a root of unity of order $p^{m+c}$, except in the special case 1 where $a=c=1$ and $r^{1-\sigma}$ has order $2^{m+l-1}$ as expected. In the next computation, we'll use 

$$m'=\left\{\begin{matrix}
 & 	m & \text{ apart from the special case}\\
   &  &  \\
 & 	m+l-2 & \text{ in the special case}\\
\end{matrix}\right.$$

It's then legitimate to write:

$$ R=r^{p^{m'+c}}= (\xi_K\delta_0^{p^{m'}})^\tau\lambda_0^{p^{m'+c}}$$

for some $\lambda_0\in\kt$.

\vx

 Now we can explicit the relationship between $N_{F/K}(t)$ and $\delta_0$ from (3.5.10): to avoid discussion about the special case 1, for an integer $k>>0$:
 
 $$N(r)^{p^k}= R^{p^{k+n-m'-c}}$$

 On the other hand, 
 
 $$N(r)=N\left(s^{\kappa S_mA_mp^{a-c}} \left(\frac{t^{S_m}}{s}\right)^{p^{n-m-c}}\right)$$
 
 $$N(r)=N(s)^{\kappa S_m(1)A_m(1)p^{a-c}-p^{n-m-c}}N(t)^{S_m(1)p^{n-m-c}}$$
 
 $$N(r)=N(s)^{\kappa A_m(1)p^{m+a-c}-p^{n-m-c}}N(t)^{p^{n-c}}$$
 
 and $A_m(1)=\sum_{k=0}^{p^{n-m}-1}k=\frac{p^{n-m}(p^{n-m}-1)}{2}$ so that 
 
  $$N(r)=N(s)^{\kappa p^{n+a-c}(p^{n-m}-1)/2-p^{n-m-c}}N(t)^{p^{n-c}}$$
  
  and using (3.5.4):
  
  $$N(r)=(\xi_Kx^{p^{n-b}})^{\kappa p^{n+a-c}(p^{n-m}-1)/2-p^{n-m-c}}N(t)^{p^{n-c}}$$

  One then has when $k>>0$, removing roots of unity:

  $$N(r)^{p^k}=R^{p^{k+n-m'-c}}=(\delta_0^{\tau p^{m'}}\lambda_0^{p^{m'+c}})^{p^{k+n-m'-c}}=(x^{p^{n-b}})^{\kappa p^{k+n+a-c}(p^{n-m}-1)/2-p^{k+n-m-c}}N(t)^{p^{k+n-c}}$$
  
  Hence,

    $$\delta_0^{\tau p^{k+n-c}}\lambda_0^{p^{k+n}}=x^{\kappa p^{2n+k+a-b-c}(p^{n-m}-1)/2-p^{2n+k-m-b-c}}N(t)^{p^{k+n-c}}$$
    
   Finally it comes:

    $$\delta_0^{\tau}\lambda_0^{p^c}=\xi_K^{-\rho} x^{\kappa p^{n+a-b}(p^{n-m}-1)/2-p^{n-m-b}}N(t)$$

     for some $\rho\in\zp$, and then:
     
     $$N(t)=\frac{\xi_K^\rho\delta_0^\tau}{x^{\kappa p^{n+a-b}(p^{n-m}-1)/2-p^{n-m-b}}} \lambda_0^{p^c}$$
     
     which can be written:
     
     $$\hspace{-5cm}(3.5.12)\hspace{3cm}N(t)=\frac{\xi_K^\rho\delta_0^\tau}{x^{\nu p^{n-m-b}}} \lambda_0^{p^c}$$
     
     for some $\nu,\ \tau\in\zp^\times$.
  
 \vx

{\bf lemma (3.5.13):}
for a convenient $x$, one has in $\kt/K^{\times p}$:

$$\dim_{\mathbb{F}_p}\langle x,\xi_K,\delta_0\rangle=3$$

and then:
 $$\dim_{\qp}(\qp\underset{\zp}{\otimes}W_0)=2$$  \vspace{0,5cm}

Proof:  
\vspace{0,5cm}

First note that in $\kt/K^{\times p}$, $\xi_K$ and $\delta_0\in K$ such that $\xi_F\delta_0=Rad(K_{m+1}/K_m) $ doesn't depend on the choice of $x$ generating $\kt\mod N(\ft)$. On the other hand, $\mu_{p^{a+1}}\not\subset K$, then $\delta_0\not\subset \mu_K K^{\times p}$, so that in $\kt/K^{\times p}$, 

$$\dim_{\mathbb{F}_p}\langle \xi_K,\delta_0\rangle=2$$

Finally, if $x$ generates $\kt\mod N(\ft)$, so is $xN(f)$ for any $f\in\ft$. Then, $\dim_{\mathbb{F}_p}\langle xN(f),\xi_K,\delta_0\rangle=2$ for any $f\in\ft$ iff both $x$ and $N(\ft)$ are contained in $\langle\xi_K,\delta_0\rangle K^{\times p}$. However, $\dim_{\mathbb{F}_p}(N(\ft))=d+1$ then one has $d=1$, $\langle\xi_K,\delta_0\rangle=N(\ft)\mod K^{\times p}$ and then $x\in N(\ft)K^{\times p}$ cannot generate $\kt\mod N(\ft)$, which is a contradiction.

Once having chosen a convenient $x$ such that $\dim_{\mathbb{F}_p}\langle x,\xi_K,\delta_0\rangle=3$, (3.5.12) proves that $$\hspace{-5cm}(3.5.14)\hspace{3cm}\dim_{\mathbb{F}_p}\langle x,\xi_K,N(t)\rangle=3$$ and then: $$\dim(\qp\underset{\zp}{\otimes}W_0)=2$$

\vspace{0,5cm}

This allows us to choose a direct factor in $N(\ft)$:
\vspace{0,5cm}

{\bf lemma (3.5.15):} a convenient direct factor \vspace{0,5cm}

$$W_0=\langle N(s),N(t),N(x)\rangle$$

is a direct factor in $N(\ft)$.\vspace{0,5cm}

Proof: \vspace{0,5cm}

At first, note that using (3.5.4), $N(s)=\xi_K x^{p^{n-b}}$ so that $N(s)^{p^b}=\xi_K^{p^b}x^{p^n}$ and then $\xi_K^{p^b}\in \langle N(s),N(x) \rangle$. On the other hand, $\xi_K^{p^b}$ generates $\mu_K\cap N(\ft) =(W_0)_{tors}$ which we denote $\mu_N$, and then $W_0$ is a direct factor in $N(\ft)$ iff $W_0/\mu_N$ is a direct factor in $N(\ft)/\mu_N$, which has no torsion. Following (1.8), we assume that, for $\xi_N=\xi_K^{p^b}$:

$$\hspace{-5cm}(3.5.16)\hspace{3cm}N(f)^p=\xi_N^\alpha N(s)^\beta x^{p^n\gamma } N(t)^\delta$$

Replacing $N(s)$ with $\xi_Kx^{p^{n-b}}$ in this equality and using (3.5.14) proves that $p\mid \delta=p\delta'$, hence we rewrite (3.5.16):
$$N(f/t^{\delta '})^p=\xi_N^\alpha N(s)^\beta x^{p^n\gamma } =\xi_K^{\alpha p^b} N(s)^\beta x^{p^n\gamma}= \xi_K^{\alpha p^b+\beta} x^{p^n\gamma +\beta p^{n-b}}$$

If $b>0$, this implies that $\xi^\beta\in K^{\times p}$ (because $n-b>0$) and then $p\mid\beta=b\beta '$, and then \newline $N(f/t^{\delta '}s^{\beta '})^p=\xi_N^\alpha x^{p^n\gamma}=\xi_K^{\alpha p^b} x^{p^n\gamma}=N(s)^{\alpha p^b}x^{(\gamma-\alpha) p^n}$, then $x^{(\gamma-\alpha) p^n}\in N(\ft)^p$, hence $p\mid (\gamma-\alpha) =pw$  so that 

$$N(f/t^{\delta '}s^{\beta '})^p=N\left(s^{\alpha p^{b-1}}x^w\right)^p$$

and then $N(f/t^{\delta '}s^{\beta '})=\xi N\left(s^{\alpha p^{b-1}}x^w\right)$
for some $\xi\in\mu_p$, necessarily norm and then in $W_0$. So we're done if $b>0$.

If $b=0$, the equality becomes 

$$N(f/t^{\delta '})^p= \xi_K^{\alpha +\beta} x^{p^n(\gamma +\beta)}$$
hence $p\mid (\beta+\alpha)=p\mu$, but since $\xi_K$ is norm (because $b=0$) it follows $\xi_K^{p\mu}\in N(\ft)^p$ and then $p\mid(\gamma+\beta)=p\nu$. But then,
$$N(f/t^{\delta '})^p= \xi_K^{p\mu} x^{p^{n+1}\nu}$$
with both factors on the right in $W_0^p$ so that we can finish as above, completing the proof of lemma (3.5.15).
\vspace{0,5cm}

The last step consists in choosing $S$ so that
$$N(\ft)=W_0\oplus S$$

with $S\simeq\zp^{d-1}$. We then can write $$S=\langle v_1,v_2,...,v_{d-1}\rangle$$
with $$v_i=N(u_i),\ 1\le i\le d-1$$
and then
$$N(\ft)=\langle N(s),N(t),N(x), N(u_1),...,N(u_{d-1})\rangle$$
so that, following the "N(L)=N(M)" lemma:
$$\ft=\langle s,t,x,u_1,...,u_{d-1}\rangle$$
Finally, let $$V=\frac{\langle X,S,T,U_1,...,U_{d-1}\rangle_{\zp[G]}}{\langle S^{1-\sigma+\kappa p^a}X^{-\kappa p^{a-b}}T^{\sigma^{p^m}-1}\rangle_{\zp[G]}}=\frac{M}{D}$$

be the formal space where, in $M$: $$\langle X,S,T,U_1,...,U_{d-1}\rangle=\langle X\rangle\oplus\langle S\rangle\oplus \langle T\rangle\oplus\langle U_1\rangle\oplus....\langle U_{d-1}\rangle  \simeq  \zp\oplus \zp[G]\oplus...\oplus\zp[G]$$

and in $V$:

$$\langle X,S,T\rangle\oplus\langle U_1\rangle\oplus....\langle U_{d-1}\rangle   \simeq    \frac{\langle X,S,T\rangle}{\langle S^{1-\sigma+\kappa p^a}X^{-\kappa p^{a-b}}T^{\sigma^{p^m}-1}\rangle}\oplus \zp[G]\oplus...\oplus\zp[G]$$
which means that $X^\sigma=X$ and $\frac{S^{1-\sigma+\kappa p^a}}{X^{-\kappa p^{a-b}}}=T^{1-\sigma^{p^m}}$ generate all relations between $X,S,T,U_1...U_{d-1}$. We study in IV the structure of $W=\frac{\langle X,S,T\rangle}{\langle S^{1-\sigma+\kappa p^a}X^{-\kappa p^{a-b}}T^{\sigma^{p^m}-1}\rangle}$, proving:
$$H^0(G,W)=p^n$$
$$H^1(G,W)=1$$
$$\chi_W=\chi_{reg}+1$$
and
$$W_{tors}=\mu_F$$
Then, $\ft$ is a quotient of $V$. But it has the same character
$$\chi_{\ft}=\chi_V=d\chi_{reg}+1$$
and has the same torsion, so that, using the usual arguments:
 $$\ft\simeq V$$

 \vspace{1cm}
 {\bf \large Case 6: } $p=2,\ a=m=1, \dim_{\mathbb{Q}_2}K $ even, $K[\mu_\infty]/K \text{ not procyclic}$, a special case \vspace{1cm}

 We dive into what happens in special case 2, assuming that $n\ge 2$ (otherwise the cyclotomic case applies). This can't happen if $\dim_{\mathbb{Q}_2}K$ is odd ( if so, -1 is not a norm in $K(i)/K$, see (3.3.5) and then $K(i)$ cannot be contained in any cyclic extension $F/K$ of degree 4). The first parameter is the number of roots of unity in $K(i)=K_1$, we then note $l\ge 2$ the integer such that:
 
 $$\mu_{K(i)}=\mu_{2^l}$$
 
 As said in III.1, $\mu_F=\mu_{K(i)}$ and we note $\xi_F$ a generator of this group. Since 
 
 $$N_{K(i)/K}(\xi_F)=\xi_F\xi_F^{-1}=1$$
 
 the commutative diagram in (1.2) applied replacing $L$ with $K(i)$ proves that $\xi_F$ is norm in any abelian extension of $K$, and so is $F$, then we can write
 
 $$\xi_F=N_{F/K_1}(z)$$
 
 for some $z\in\ft$. Moreover, 
 
 $$N(z)=N_{K_1/K}(\xi_F)=1$$
 
 so that 
 
 $$z=t^{1-\sigma}$$
 
 for some $t\in\ft$. On the other hand,
 
 $$1=\xi_F^{2^l}=N_1(t^{2^l(1-\sigma)})$$
 
 where $N_1=N_{F/K_1}=\underset{\hspace{-0,3cm}0\le k \le 2^n-1}{\sum \ \sigma^{2k}} $ and applying H90 in $F/K_1$:
 
 $$\hspace{-5cm}(3.6.1)\hspace{3cm}t^{2^l(1-\sigma)}=y^{1-\sigma^2}$$
 
 for some $y\in\ft$. Then, one has
 
 $$\left(\frac{t^{2^l}}{y^{1+\sigma}}\right)^{1-\sigma}=1$$
 
 and then $$\frac{t^{2^l}}{y^{1+\sigma}}=x_0\in\kt$$
 
 We then notice that $y$ found in (3.6.1) can be replaced with $y\delta_1$ for any $\delta_1\in K_1$ since $(y\delta_1)^{1-\sigma^2}=y^{1-\sigma^2}$. Then, $x_0$ becomes
 
 $$x'_0=x_0/N_{K_1/K}(\delta_1)$$
 
 and then $x_0$ can be changed modulo $N_{K_1/K}(K_1^\times)$, of which the index in $\kt$ is 2. This will be usefull later.
 
 $$\hspace{-5cm}(3.6.2)\hspace{3cm}t^{2^l}=x_0 y^{1+\sigma}$$
 
 is our fundamental relation, to be more precise later.
 
 The next step consists in determining what $N(t)$ and $N(y)$ may look like. We start with
 
 $$N_1(z)=\xi_F=N_1(t)^{1-\sigma}$$
 
 and then notice that
 
 $$(1-\xi_F)^{1-\sigma}=\frac{1-\xi_F}{1-\xi_F^{-1}}=-\xi_F$$
 
 Comparing the 2 previous equalities leads to:
 
 $$N_1(-it)^{1-\sigma}=-\xi_F=(1-\xi_F)^{1-\sigma}$$
 
 and then,
 
 $$\hspace{-5cm}(3.6.3)\hspace{3cm}N_1(t)=(1-\xi_F)i\delta_0$$
 
 for some $\delta_0\in\kt$. Then,
 
 $$N(t)=((1-\xi_F)i\delta_0)^{1+\sigma}$$
 
 $$\hspace{-5cm}(3.6.4)\hspace{3cm}N(t)=(2-2\cos\theta)\delta_0^2$$
 
 where we note $\theta=\frac{2\pi}{2^l}$ and $N_{K_1/K}(1-\xi_F)=2-2\cos\theta\in\kt$. Of course, a finite extension of $\mathbb{Q}_2$ is not the same as the complex plane, but the following calculation is still legitimate (think about both $F$ and $\mathbb{C}$ as completions of the same number field):
 
 $$1-e^{i\theta}=e^{i\theta/2}(e^{-i\theta/2}-e^{i\theta/2})=e^{i\theta/2}(-2i\sin(\theta/2)$$
 
 and then:
 
 $$(1-e^{i\theta})^2=-4e^{i\theta}\sin^2(\theta/2)=-e^{i\theta}(2-2\cos\theta)$$
 
 Notably:
 
 $$\hspace{-5cm}(3.6.5)\hspace{3cm}(1-\xi_F)^{2^l}=-(2-2\cos\theta)^{2^{l-1}}$$
 
 Now, applying $N_1$ in (3.6.2):
 
 $$N_1(t)^{2^l}=x_0^{2^{n-1}}N(y)$$
 
 along with (3.6.3) and (3.6.5):
 
 $$((1-\xi_F)i\delta_0)^{2^l}=x_0^{2^{n-1}}N(y)$$
 
 and then:

 $$\hspace{-5cm}(3.6.6)\hspace{3cm}N(y)=\frac{-(2-2\cos\theta)^{2^{l-1}}}{x_0^{2^{n-1}}}\delta_0^{2^l}$$
 
 Once again from (3.6.2), taking norms:
 
 $$N(t)^{2^l}=x_0^{2^n}N(y)^2$$
 
 and then,
 
 $$N(t)^{2^{l-1}}=\pm \ x_0^{2^{n-1}}N(y)$$
 
 Determining the sign in this equality is now easy: if it was a plus, $N(y)$ would be a square in $K$ (because $n\ge 2$) and then using (3.6.6), $-1$ would be a square in $K$, that is $i\in\kt$. This is a contradiction, so we can write:
 
 $$\hspace{-5cm}(3.6.7)\hspace{3cm}N\left(\frac{t^{2^{l-1}}}{y}\right)=-\ x_0^{2^{n-1}}$$
 
 Remembering that $x_0$ is modulo $N_{K_1/K}(\kt_1)$, either $x_0\in N_{K_1/K}(\kt_1)$ and we'll go on with $x_0=1$; either $x_0\notin N_{K_1/K}(\kt_1)$ and we'll go on with $x_0$ being any generator of $\kt/N_{K_1/K}(\kt_1)$ and then, as well a generator of $\kt/N(\ft)$. In the first case, (3.6.7) implies that $-1$ is a norm in $F/K$; in the second, $-1$ is not a norm, because otherwise $x_0^{2^{n-1}}$ would be norm as well and couldn't generate $\kt/N(\ft)$. (3.6.2) then becomes:
 
  $$\hspace{-1cm}(3.6.8)\hspace{3cm}\left\{ \begin{matrix}
 & t^{2^l}=y^{1+\sigma}& \text{ if } -1 \text{ is a norm in}\ F/K \\
 &  & \\
 & t^{2^l}=x_0 y^{1+\sigma}& \text{ if } -1 \text{ is not a norm in }\ F/K
  \end{matrix} \right.$$

Note that in the second case, $x_0$ denotes a generator of $\kt/N(\ft)$. We then need the following lemma to find a convenient direct factor in $N(\ft)$:

\vx

{\bf lemma (3.6.9):} None of $-1,\ 2-2\cos\theta, -2+2\cos\theta$ is a square in $\kt$.

\vx

Proof: We already know that $-1$ is not a square, and can use the formula $2-2\cos\theta=4\sin^2(\theta/2)$ to see that if $2-2\cos\theta$ or $-2+2\cos\theta$ was a square, $K(i)$ would contain $i$, $\sin(\theta/2)$ and $\sin\theta=2\sin(\theta/2)\cos(\theta/2)$ so that $e^{i\theta/2}$, a root of unity of order $2^{l+1}$ would be in $F$, which is a contradiction.

\vx

It follows from (3.6.4) and (3.6.6) that in $\kt/K^{\times 2}$:

$$\langle N(t), N(y) \rangle_{\mathbb{F}_2} = \langle -1, 2-2\cos\theta\rangle$$

and then, according to the previous lemma, 

$$W_0=\langle N(t), N(y) \rangle_{\mathbb{Z}_2}$$

is a direct factor in $N(\ft)$. In the first case where $-1$ is a norm, one has $-1=N(t^{2^{l-1}}/y)$. In the second case, $N(\ft)$ has no torsion; so that 

 $$\left\{ \begin{matrix}
	& \dim\mathbb{Q}_2\underset{\mathbb{Z}_2}{\otimes}W_0=1& \text{ if } -1 \text{ is a norm in}\ F/K \\
	&  & \\
	&\dim\mathbb{Q}_2\underset{\mathbb{Z}_2}{\otimes}W_0=2& \text{ if } -1 \text{ is not a norm in }\ F/K
\end{matrix} \right.$$

It's easy to see that in the first case where $-1$ is a norm, $x_0$ being any generator of $\kt/N(\ft)$, $\langle x_0,N(t), N(y) \rangle_{\zp}$ is as well a direct factor in $N(\ft)$: basically, if $N(f^2)=\pm\  x^{2^n\alpha}N(t)^\beta N(y)^\gamma$, we can suppose $\gamma=0$ using $-1=N(t^{2^{l-1}}/y)$ and lemma (3.6.9) proves that $2\mid\beta$; it follows $N(h^2)=\pm\ x^{2^n\alpha}$ for some $h\in\ft$, then it's a plus (because $-1$ is not a square) and the reader can conclude that $2\mid\alpha$. 

\vx

This leads to:

 $$\left\{ \begin{matrix}
	& V=\langle X\rangle\oplus\langle T,Y\rangle\oplus\langle U_1\rangle\oplus....\langle U_{d-1}\rangle   \simeq   \mathbb{Z}_2\oplus \frac{\langle T,Y\rangle}{\langle T^{2^l}Y^{-1-\sigma}\rangle}\oplus \mathbb{Z}_2[G]\oplus...\oplus\mathbb{Z}_2[G]& \text{ if } -1 \text{ is a norm in}\ F/K \\
	&  & \\
	& V=\langle X,T,Y\rangle\oplus\langle U_1\rangle\oplus....\langle U_{d-1}\rangle   \simeq    \frac{\langle X,T,Y\rangle}{\langle T^{2^l}X^{-1}Y^{-1-\sigma}\rangle}\oplus \mathbb{Z}_2[G]\oplus...\oplus\mathbb{Z}_2[G]& \text{ if } -1 \text{ is not a norm in }\ F/K
\end{matrix} \right.$$

 Using the usual arguments, after studying the formal spaces
 
 $$W_1=\frac{\langle X,T,Y\rangle}{\langle T^{2^l}Y^{-1-\sigma}\rangle}$$
 
 and
 
  $$W_2=\frac{\langle X,T,Y\rangle}{\langle T^{2^l}X^{-1}Y^{-1-\sigma}\rangle}$$
  
  We'll be able to write:

$$\ft\simeq V$$

\vspace{1cm}
{\bf \large Case 7: } $Char(k)\neq p$
 \vspace{1cm}

This last case considers an extension $F/K$ of degree $p^n$ when the residual field $k_K$ of $K$ has $q=q_0^{i}$ elements, $q_0$ being a prime number with $q_0\neq p$. In this case, one has (since we consider $p$-parts only:)

$$\hspace{-5cm}(3.7.1)\hspace{3cm}\kt=\mu_K\times \pi_K^{\zp}$$

for any uniformizer $\pi_K$ of $K$.

\vx

Either $\mu_K=1$, then $F/K$ is not ramified (because $K^{ab}$, the maximal pro-p extension abelian over $K$, is unramified) and $p\neq 2$; moreover, the number of automorphisms of $\mathbb{Z}/p\mathbb{Z}$ is prime to $p\neq 2$, so that $\mu_F=1$ and then, remembering that we only consider $p$-parts:

$$\ft=\kt$$
\vx

Either $\mu_K\neq 1$. Since $\mu_K$ corresponds to the ramification subgroup of $\gal(K^{ab}/K)$ in Class Field Theory, 

$$p^b=e$$

is the ramification index of $F/K$, while

$$p^m=f$$

is the inertia index. We then have 

$$n=m+b$$

We then notice that $c=\min(n-m,a)=\min(b,a)=b$, then $K_{m+c}=K_n=F$. Lemma (3.5.8) and remark (3.5.9) apply and prove that $F/K$ is pseudo-Kummer in the sense that 

$$F=K\left[\sqrt[p^n]{r}\right]$$

for some $r\in\kt\backslash K^{\times p}$. $K_m/K$ is not ramified, then (maybe replacing $r$ with $r^\lambda$ with $\lambda\in\zp^\times$), we can assume that 

$$r=\xi_K^\alpha\pi_K^{p^m}$$

for some $\alpha\in\zp^\times$ and some uniformizer $\pi_K$ of $K$. At last, we note 

$$\hspace{-5cm}(3.7.2)\hspace{3cm}\pi_F=\sqrt[p^n]{\xi_K^\alpha\pi_K^{p^m}}$$

which is a uniformizer of $F$. Except when $p=2$ and $a=1$, one has 
$$|\mu_F|=p^{a+m}$$

 If $p=2$ and $a=1$, remembering that the Frobenius automorphism sending $x\in k_K$ to $x^q$ generates $\gal(K_m/K)=\gal(k_F/k_K)$, one has $2\mid (q-1)$ and $4\nmid (q-1)$. Then, $q=-1+\kappa 2^{l-1}$ for some $\kappa, l\in \mathbb{N}$ with $\kappa$ odd and $l\ge 3$ (because $q\equiv -1\ [4]$). It follows $q^2=1-\kappa 2^{l}+\kappa^2 2^{2l-2}$ so that if $m=1$,  $\mu_F=\mu_{2^{l}}$ and an easy induction proves 
 
 $$|\mu_F|=2^{l+m-1}$$
 
 like in special case 1. Anyway, in all cases, according to (3.7.1) in $F$:

$$\hspace{-5cm}(3.7.3)\hspace{3cm}\ft=\langle \pi_F,\mu_F\rangle_{\zp}$$

and 

$$\hspace{-5cm}(3.7.4)\hspace{3cm}\pi_F^{1-\sigma}=\xi$$

a root of unity of order $p^n$ (according to (3.7.2)). We then define $\delta$ the positive integer such that 

$$p^\delta = \left\vert\frac{\mu_F}{\langle\xi\rangle }\right\vert=\left\{\begin{matrix}
 &	p^{a-b} & \text{ if $p\neq 2$ or $a\neq 1$}\\
&	 & \\
&	 2^{l-b-1} & \text{ if $p= 2$ and $a= 1$}
\end{matrix}\right.$$

We choose $\xi_F$ generating $\mu_F$ such that:

$$\hspace{-5cm}(3.7.5)\hspace{3cm}\xi=\xi_F^{p^\delta}$$

and at last, we note $k_\sigma\in\zp^\times$ a p-adic unity such that:

$$\forall \zeta\in\mu_F,\  \zeta^\sigma=\zeta^{k_\sigma}$$

Then, using (3.7.4), the action of $G$ on $\ft$ is then given by:

$$\hspace{-5cm}(3.7.6)\hspace{3cm}\left(\pi_F^\alpha\xi_F^\beta\right)^\sigma=\pi_F^\alpha\xi_F^{p^\delta \alpha +k_\sigma\beta}$$

for all $\alpha$ and $\beta$; choosing them in $\zp$ is enough to completely describe the action of $G$ according to (3.7.3), but they can as well be chosen in $\zp[G]$ (because its commutative).

\vspace{1cm}

{\bf \large IV Torsion submodules and Cohomology of Formal Spaces } \vspace{1cm}

In this section, we establish basic properties of the formal spaces introduced in III. We begin with the most general space.

Let $G=\mathbb{Z}/p^n\mathbb{Z}$ be the cyclic group of order $p^n$ and $W_{a,b,m,n}$ be the formal space defined as a $\zp[G]$-module for any positive integers $a,b,m,n$ satisfying:
\vx

$$\left\{\begin{matrix}
	 & n\ge 1 \\ & \\
	  & a\ge 1 \\ & \\
	 & b\le\min(a,n)\\ & \\
 &	m+b\le n
\end{matrix}
\right.$$

 by:

$$W_{a,b,m,n}=\frac{\langle X,S,T\rangle_{\zp[G]}}{\langle S^{-\sigma+1+\kappa p^a}X^{-\kappa p^{a-b}}T^{\sigma^{p^m}-1}\rangle_{\zp[G]}}=\frac{M}{D}$$

\vx

where $X^\sigma=X$ and $S^{-\sigma+1+\kappa p^a}=X^{\kappa p^{a-b}}T^{1-\sigma^{p^m}}$ generate all relations between $X,S,T$ (the first relation holds both in $M$ and $W=W_{a,b,m,n}$ while the second holds in $W$ only) and where $\kappa\in\zp^\times$ with the following restriction: if $p=2$ and $a=1$, we must have:

$$\kappa\neq -1$$

otherwise calculation of $\chi(W)$ would be wrong. The reader can check that in special case 1, $\kappa\equiv 1 \ [4]$ and then the previous condition is satisfied, while in special case 2, $\kappa\equiv -1 \mod |\mu_F|/2$ and then, another formal space will be used as seen in Case 6. 

The reader should note that, even if $W_{a,b,m,n}$ is especially dedicated to the Mixed Case, we'll see that we can apply it in some more cases, and then most restrictions of Case 5 ($m>0$, $n>m$) are not considered here.
\vx

\begin{enumerate}
	\item {\bf Character of $W$:}\vx
	
	To determine $\chi(W)$, we note that $P\in\qp[G]$ has an inverse in $\qp[G]$ iff:
	
	$$\forall\xi\in\mu_{p^n},\ P(\xi)\neq 0$$
	
	where $P(\xi)$ is obtained replacing $\sigma$ with $\xi$ in the expression of $P$, which makes sense for any $\xi\in\mu_{p^n}$. We then note 
	
	$$q=S^{-\sigma+1+\kappa p^a}X^{-\kappa p^{a-b}}T^{\sigma^{p^m}-1}$$
	
	a generator of $D$, so that 
	
	$$\chi(W)=\chi(M)-\chi (D)=2\chi_{reg}+1-\chi(D)$$
	
	because $$M\simeq \zp\oplus\zp[G]\oplus\zp[G]$$
	
	We now assume that $q^\alpha=1$ in $\qp\underset{\zp}{\otimes}M$ with $\alpha\in\qp[G]$. Looking at the $S$-part, one has 
	
	$$(-\sigma+1+\kappa p^a)\alpha =0$$
	
	in $\qp[G]$, so that for all $\xi\in\mu_{p^n}$:
	
	$$(-\xi+1+\kappa p^a)\alpha(\xi) =0$$
	
	in $\qp$. The only root of $-X+1+\kappa p^a$ is $r=1+\kappa p^a$. If $p\neq 2$, $r\in\zp$ hence cannot be a root of unity of order a power of $p$. If $p=2$, $\mu_{2^n}\cap\mathbb{Z}_2=\mu_2$ and $r=1$ is excluded (because $\kappa\in\mathbb{Z}_2^\times$) while $r=-1$ is possible iff $\kappa\times 2^a=-2$ that is $\kappa=-1$ and $a=1$. We excluded this case, and then finally, $-\sigma+1+\kappa p^a$ is invertible in $\qp[G]$, hence:
	
	$$\forall\xi\in\mu_{p^n},\ \alpha(\xi)=0$$
	
	and then $\alpha=0$ in $\qp[G]$, so that $\chi(D)=\chi_{reg}$ and one has:
	
	$$\hspace{-5cm}(4.1.1)\hspace{3cm}\chi(W)=\chi_{reg}+1$$\vx

			\item {\bf more general results}\vx		
			
			We make computations, which we will allow to compute $H^i(G,W)$ for $i=0,1$ and also $W_{tors}$, both for the most general formal space $W_{a,b,m,n}$ and for the formal spaces $W_1,\ W_2$ as defined in Case 6. \vx
			
			We define

		$$W=\frac{\langle X,S,T\rangle_{\zp[G]}}{\langle S^{P}X^{Q}T^R\rangle_{\zp[G]}}=\frac{M}{D}$$
		
		for some $P,Q,R\in\zp[G]$, where $X^\sigma=X$ and $S^PX^QT^R=1$ generate all relations between $X,S,T$ (the first relation holds both in $M$ and $W$ while the second holds in $W$ only). We note $q=S^PX^QT^R\in M$.
		
		We assume moreover that:
		
		$$\forall\xi\in\mu_{p^n},\ P(\xi)\neq 0$$
		
		which notably implies (like in IV.1) that 
		$$\chi(W)=2\chi_{reg}+1-\chi_{reg}=\chi_{reg}+1$$
		
		We then consider $B$ and $C$ such that, in $\zp[X]$:
		
		$$BC=X^{p^n}-1$$
		
		which notably implies that $B$ and $C$ are primitive and at last, we define:
		
		$$H(B,C,W)=\frac{Ker\  (B,W)}{Im\  (C,W)}=\frac{\{w\in W,\ w^B=1\}}{\{w\in W,\ \exists v\in W,\ w=v^C\}}$$
		
		and we compute $H(B,C,W)$ this way: assuming $w\in Ker\ (B,W)$, one has in $M$:
		
		$$w^B=q^\lambda$$
		
		for some $\lambda\in\zp[G]$ and then if $w=S^\alpha X^\beta T^\gamma$ for some $\alpha, \beta, \gamma\in\zp[G]$:
		
		$$w^B=S^{\alpha B} X^{\beta B} T^{\gamma B}=S^{P\lambda}X^{Q\lambda}T^{R\lambda}$$
		
		On the $S$ side, it follows:
		
		$$\forall\xi\in\mu_{p^n},\ B(\xi)=0\implies P(\xi)\lambda(\xi)=0\implies \lambda(\xi)=0$$
		
		because $P$ was supposed to have no root in $\mu_{p^n}$. But then $B\mid \lambda=B\lambda '$ and:
		
		$$w^B=q^{B\lambda '}$$
		
		so that 
		
		$$w=q^{\lambda '}\delta$$
		
		for some $\delta\in M$ such that $\delta^B=1$, that is $\delta\in Ker (B,M)$. One has:
		
		$$Ker (B,M)=\left\{\begin{matrix}
		 & 	S^{C\zp[G]}X^{\zp[G]}T^{C\zp[G]} & \text{ if }B(1)=0\\
		  & & \\
		 
		 	 & 	S^{C\zp[G]}T^{C\zp[G]} & \text{ if }B(1)\neq0
		\end{matrix} \right.$$
	
	and then, since $w=\delta$ in $W$ with $\delta\in Ker(B,M)$:
	
	$$H(B,C,W)=\left\{\begin{matrix}
		& 	\frac{\langle X\rangle}{\langle X\rangle\cap W^C} & \text{ if }B(1)=0\\
		& & \\
		
		& 	1 & \text{ if }B(1)\neq0
	\end{matrix} \right.$$
	
	Computing $\langle X\rangle\cap W^C$ is quite easy: if
	
	$$X^\mu=S^{\alpha C} X^{\beta C} T^{\gamma C}S^{P\lambda}X^{Q\lambda}T^{R\lambda}$$
	
	the $S$ side proves as above that $C\mid\lambda$ and then $X^\mu\in\langle X\rangle^C$. Finally, it comes:
	
		$$H(B,C,W)=\left\{\begin{matrix}
		& 	\frac{\langle X\rangle}{\langle X\rangle^{C(1)}} & \text{ if }B(1)=0\\
		& & \\
		
		& 	1 & \text{ if }B(1)\neq0
	\end{matrix} \right.$$

of which the order is:

	$$\hspace{-3cm}(4.2.1)\hspace{3cm}\vert H(B,C,W)\vert =\left\{\begin{matrix}
	& \vert\ 	\zp / C(1)\zp\ \vert & \text{ if }B(1)=0\\
	& & \\
	
	& 	1 & \text{ if }B(1)\neq0
\end{matrix} \right.$$

Note that taking $B=1-\sigma,\ C=N$ and $B=N,\ C=1-\sigma$ leads to calculation of $H^0(G,W)$ and $H^1(G,W)$ respectively, so that

	$$\hspace{-5cm}(4.2.2)\hspace{3cm}\vert H^i(G,W)\vert =\left\{\begin{matrix}
	& p^n & \text{ if }\ i=0\\
	& & \\
	
	& 	1 & \text{ if }\ i=1
\end{matrix} \right.$$\vx

Now, we aim to study $W_{tors}$. To do so, we'll make another assumption about $R$: from now on, in $\zp[X]$:

$$X^{p^n}-1=RV$$

with both $R$ and $V$ primitive polynomials. We then assume that $w=S^\alpha X^\beta T^\gamma \in W_{tors}$, that is for some $\lambda\in\zp$:

$$\hspace{-3cm}(4.2.3)\hspace{3cm}w^{p^k}=S^{p^k\alpha}X^{p^k\beta}T^{p^k\gamma}=q^\lambda=S^{P\lambda}X^{Q\lambda}T^{R\lambda}$$

We make the euclidean division of $\lambda$ by $V$:

$$\lambda=Vt+\rho$$

with $d^\circ t<d^\circ V$. Since $RV=0$ in $\zp[G]$, the $T$ side of (4.2.3) proves that:

$$p^k\gamma = R\rho$$

with $R$ primitive and $d^\circ (R\rho)<d^\circ R+d^\circ V=p^n$, so that:

$$\rho=p^k\rho'$$

in $\zp[G]$, then $\lambda=Vt+p^k\rho'$  and (4.2.3) becomes:

$$(w/q^{\rho '})^{p^k}=S^{PVt}X^{QVt}=q^{Vt}$$

which implies that $p^k\mid PVt$ and $p^k\mid QVt$ and then for some $\mu, \nu\in\zp[G]$:

$$w\equiv S^{\mu V}X^{\nu}\ \mod q$$

We now determine the conditions for $w=S^{\mu V}X^{\nu}\in W$ to be in $W_{tors}$. At first, $P$ is invertible in $\qp[G]$ and then, there exists $C\in\zp[G]$ such that $PC=p^j$ in $\qp[G]$. Hence, $$w\in W_{tors}\iff w^P\in W_{tors}$$

Now, $w^P=S^{\mu VP}X^{\nu P}$ is a torsion element iff for $i>>0$, in $M$:

$$w^{p^iP}=S^{p^i\mu VP}X^{p^i\nu P}=S^{P\lambda}X^{Q\lambda}$$

for some $\lambda\in\zp[G]$. The $S$ side along with the fact that $P$ is invertible in $\qp[G]$ proves that $\lambda=p^i\mu V$ and the $X$ side then provides (remember that $P(1)\neq 0$) in $\qp$:

$$\nu(1)=\frac{Q(1)V(1)}{P(1)}\mu(1)$$

If this equality holds, one has $w^{p^iP}=q^{p^i\mu V} $, so that finally:

$$w=S^{\mu V}X^{\nu }\in W_{tors}\iff \nu(1)=\frac{Q(1)V(1)}{P(1)}\mu(1)$$

This proves that $W_{tors}$ is generated by:

$$\hspace{-2cm}(4.2.4)\hspace{2cm}\begin{matrix}
	&\Xi=S^V&\text{ if V(1)Q(1)=0}\\
	 & & \\
	& \Xi=S^VX^{\frac{Q(1)V(1)}{P(1)}} &\text{ if $P(1)\mid Q(1)V(1)\neq 0$} \\
	& & \\
	& \Xi=S^{\frac{P(1)}{Q(1)V(1)}V}X\text{ and }\Xi'=S^{V(1-\sigma)} &\text{ if $0\neq Q(1)V(1)\mid P(1)$}
\end{matrix}$$

At last, in $W$, one has $\Xi^P=\Xi^R=1$. Conversely, $\Xi^\alpha=1$ in $W$ implies on the $S$ side that $P\mid V\alpha$ in $\zp[G]$. We can consider that $\alpha$ is modulo $R$, because $S^{RV}=1$; but then, assuming $d^\circ\alpha<d^\circ R$, one has $d^\circ(V\alpha)<p^n$ and then $P\mid V\alpha$ in $\zp[X]$, and then $P\mid\alpha$ (because $P$ has no roots in $\mu_{p^n}$ and $V\mid (X^{p^n}-1)$). Then, the annihilator ideal of $\Xi$ is:

$$I=(P,R)$$

and in the third case of (4.2.4), the annihilator of $\Xi'$ is, observing that in that case $RV=X^{p^n}-1$ and $Q(1)V(1)\neq 0$ implies that $(X-1)\mid R=(X-1)R'$:

$$J=(P,R')$$

It follows:

$$\hspace{-2cm}(4.2.5)\hspace{2cm}W_{tors}\simeq\left\{ \begin{matrix}
& \zp[X]/(P,R)&\text{if $P(1)\mid V(1)Q(1)$ ($=0$ or not)}\\
& & \\
& (\frac{P(1)}{Q(1)V(1)},1-X)/(P,R)&\text{if $0\neq Q(1)V(1)\mid P(1)$}
\end{matrix}\right.$$

The isomorphism simply sends $S^{\alpha V}$ (eventually multiplied by a power of $X$) to $\alpha$. Note that in the second case of (4.2.5), the ideal in the numerator contains the one in the denominator because $(1-X)\mid R$ and $P=P(1)+(1-X)Z= \frac{P(1)}{Q(1)V(1)}\times Q(1)V(1)+(1-X)Z$ for some $Z\in\zp[X]$.		
		
				\vx

\item {\bf Application to $W_{a,b,m,n}$:}\vx				
	
	$$W_{a,b,m,n}=\frac{\langle X,S,T\rangle_{\zp[G]}}{\langle S^{-\sigma+1+\kappa p^a}X^{-\kappa p^{a-b}}T^{\sigma^{p^m}-1}\rangle_{\zp[G]}}$$
	
	And then we apply here the results found in IV.2 with:
	
	$$\begin{matrix}
		P=-X+1+\kappa p^a\\
		Q=-\kappa p^{a-b}\\
		R=X^{p^m}-1 \\
		V=N_m=\sum_{k=0}^{p^{n-m}-1}X^{kp^m}
	\end{matrix}$$

(4.2.2) then proves that:

$$|H^0(G,W_{a,b,m,n})|=p^n$$

and

$$|H^1(G,W_{a,b,m,n})|=1$$	

as expected.\vx

To describe $W_{tors}$, we first note that $P(1)=\kappa p^a$, $Q(1)=-\kappa p^{a-b}$ and $V(1)=p^{n-m}$, so that $P(1)\mid Q(1)V(1)$ because $n-m-b\ge 0$ (see (3.5.2)). Then, (4.2.4) states that the torsion submodule of $W_{a,b,m,n}$ is generated by

$$\Xi=S^{N_m}X^{-p^{n-m-b}}$$

which was supposed to be equal to $\xi_F$ in $\ft$; and (4.2.5) gives:

$$W_{tors}\simeq \zp[X]/(X^{p^m}-1,-X+1+\kappa p^a)$$

	We now have to compare this with $\mu_F$, which is $\omega$-isotypical as well. Special case 2 is excluded, then 

$$\mu_F=\mu_{K_m}=\left\{
\xi\in\mu_\infty, \xi^{\sigma^{p^m}}=\xi\right\}$$

(in special case 2, the set on the right may be larger than $\mu_F$) and then, $$\sigma^{p^m}-1=-\sigma+1+\kappa p^a=1$$ holds as well in $\mu_F$ which is then a quotient of $\zp[X]/(X^{p^m}-1,-X+1+\kappa p^a)$. The reader can conclude that this quotient is as well an isomorphism, checking that the order on both sides is $p^{a+m}$ except in special case 1 where it is $2^{m+l-1}$ (and not $2^{m+1}$). Finally,

$$\hspace{-3cm}(4.3.1)\hspace{3cm}W_{tors}\simeq \mu_F\simeq \zp[X]/(X^{p^m}-1,-X+1+\kappa p^a)$$

as expected.

				\vx

		\item {\bf Discussion about various cases:}\vx		
		
		Previous computations proves that the isomorphism 
		
		$$\ft\simeq V$$
		
		found in Case 5 is correct. What about the other ones ? \vx
		
		\begin{enumerate}
			\item \underline{ $m=0,\ \mu_F=\mu_K$:}
			
			\vx
			
			$$W_{a,b,0,n}=\frac{\langle X,S,T\rangle_{\zp[G]}}{\langle S^{-\sigma+1+\kappa p^a}X^{-\kappa p^{a-b}}T^{\sigma-1}\rangle_{\zp[G]}}$$
			
			Let's define $Z=S^\kappa$, $X_0=X^\kappa$ and $Y=T/S$: one has $\langle X,S,T\rangle=\langle X_0,Z,Y\rangle$ and
			
			$$S^{-\sigma+1+\kappa p^a}X^{-\kappa p^{a-b}}T^{\sigma-1}=Z^{p^a}X_0^{-p^{a-b}}Y^{\sigma-1}$$
			
			and then:
			
			$$W_{a,b,0,n}=\frac{\langle X_0,Y,Z\rangle_{\zp[G]}}{\langle Z^{p^a}X_0^{-p^{a-b}}Y^{\sigma-1}\rangle_{\zp[G]}}$$

			which is exactly the formal space to study in Case 4, and so it is: the isomorphism is proved as well in that case. \vx
			
			\item \underline{ $b=0,\ \mu_K\subset N(\ft)$:}
			
$$W_{a,0,m,n}=\frac{\langle X,S,T\rangle_{\zp[G]}}{\langle S^{-\sigma+1+\kappa p^a}X^{-\kappa p^{a}}T^{\sigma^{p^m}-1}\rangle_{\zp[G]}}$$

Let's define $Z=S/X$: 

$$S^{-\sigma+1+\kappa p^a}X^{-\kappa p^{a}}T^{\sigma^{p^m}-1}=Z^{-\sigma+1+\kappa p^a}T^{\sigma^{p^m}-1}$$		

and then, 

$$W_{a,0,m,n}=\frac{\langle X,Z,T\rangle_{\zp[G]}}{\langle Z^{-\sigma+1+\kappa p^a}T^{\sigma^{p^m}-1}\rangle_{\zp[G]}}=\langle X\rangle\oplus \frac{\langle Z,T\rangle_{\zp[G]}}{\langle Z^{-\sigma+1+\kappa p^a}T^{\sigma^{p^m}-1}\rangle_{\zp[G]}}$$

This is the most interesting situation, because $W$ splits into a direct summand. It doesn't correspond to any Case studied in III, but had to be mentionned. This happens as soon as $\mu_K\subset N(\ft)$, notably in $\zp$-extensions.\vx

	\item \underline{ $m=0$ and $b=0,\ \mu_K=\mu_F\subset N(\ft)$:}

	$$W_{a,0,0,n}=\frac{\langle X,S,T\rangle_{\zp[G]}}{\langle S^{-\sigma+1+\kappa p^a}X^{-\kappa p^{a}}T^{\sigma-1}\rangle_{\zp[G]}}$$
	
	Define $Z=S^\kappa/X^\kappa$, $Y=T/S$, $X_0=X$:
	
	$$W_{a,0,0,n}=\frac{\langle X_0,Y,Z\rangle_{\zp[G]}}{\langle Z^{p^a}Y^{\sigma-1}\rangle_{\zp[G]}}$$
	
	and this is the space $W$ (which splits as seen in (b)) to be studied in Case 2, of which the isomorphism is then validated. \vx
	
		\item \underline{ $m=n$ and $b=0,\ \mu_K\subset N(\ft)$ and $F/K$ is cyclotomic:}
		
		$$W_{a,0,n,n}=\frac{\langle X,S,T\rangle_{\zp[G]}}{\langle S^{-\sigma+1+\kappa p^a}X^{-\kappa p^{a}}\rangle_{\zp[G]}}$$
		
		Define $Z=S^\kappa/X^\kappa$:
		
		$$W_{a,0,n,n}=\frac{\langle X,Z,T\rangle_{\zp[G]}}{\langle Z^{-\sigma+1+\kappa p^a}\rangle_{\zp[G]}}\simeq \langle X\rangle\oplus\frac{\langle Z\rangle}{\langle Z^{-\sigma+1+\kappa p^a}\rangle}\oplus\langle T\rangle\simeq \zp\oplus\omega\oplus\zp[G]$$
		
		There was no formal space to study in the cyclotomic case, however we can note that $W_{a,0,n,n}$ as described above is really a direct factor of $\ft$ as described in Case 3, except in special case 2.\vx

		\item \underline{ $a=b=m=0,\ \mu_K=1$:}

	$$W_{0,0,0,n}=\frac{\langle X,S,T\rangle_{\zp[G]}}{\langle S^{-\sigma+1+\kappa }X^{-\kappa }T^{\sigma-1}\rangle_{\zp[G]}}$$
	
	Define $Z=S^\kappa/X^\kappa$ and $Y=T/S$:
	
	$$W_{0,0,0,n}=\frac{\langle X,Y,Z\rangle_{\zp[G]}}{\langle ZY^{\sigma-1}\rangle_{\zp[G]}}=\langle X,Y\rangle\simeq \zp\oplus\zp[G]$$
	
	Once again there was no space to study in Case 1 but $W_{0,0,0,n}$ is a direct factor as well.

	\vx

	\item \underline{ a special case: $p=2,\ a=1,\ K[\mu_\infty]/K$ not procyclic}
	\vx
	
	This is special case 2 studied in Case 6, and we have to study
	
 $$W_1=\frac{\langle X,T,Y\rangle}{\langle T^{2^l}Y^{-1-\sigma}\rangle}$$

and

$$W_2=\frac{\langle X,T,Y\rangle}{\langle T^{2^l}X^{-1}Y^{-1-\sigma}\rangle}$$

At first, computation of $H^i(G,W_j)$ for $i,j\in\{1,2\}$ is given by (4.2.2) and matches with expectations. 

After that, we can compare $W_1$ and $W_2$ with $W_{l,b,0,n}$ obtained for $p=2$ (we changed the name of one variable):

$$W_{l,b,0,n}=\frac{\langle X,S,U\rangle_{\mathbb{Z}_2[G]}}{\langle S^{-\sigma+1+\kappa 2^l}X^{-\kappa 2^{l-b}}U^{\sigma-1}\rangle_{\mathbb{Z}_2[G]}}$$

We can observe than with the following change of variables: $Y=S/U$, $T=X^\kappa/S^\kappa$, one has, replacing $b$ with 0 and 1 respectively:

$$W'_1\simeq W_{l,0,0,n}$$

and 
$$W'_2\simeq W_{l,1,0,n}$$

where 

 $$W'_1=\frac{\langle X,T,Y\rangle}{\langle T^{2^l}Y^{-1+\sigma}\rangle}$$

and

$$W'_2=\frac{\langle X,T,Y\rangle}{\langle T^{2^l}X^{-1}Y^{-1+\sigma}\rangle}$$

To summarize, $W_1$ and $W_2$ are quite similar to $W_{l,b,0,n}$ with $b=0,1$ respectively, except that the sign above $\sigma$ is unfortunately not the same in the denominator, which of course changes a lot of things. We then give up the comparison, even if $W_1$ can be compared, after isolating the $X$-part which splits, with a tensorized version of $W_{l,0,0,n}$.

To compute $W_{i,tors}$ for $i=1,2$ we then apply (4.2.5) with $P=2^l$, $Q=-b$ ($Q=0$ for $W_1$, $Q=-1$ for $W_2$), $R=-1-X$, $V=\sum_{k=0}^{2^n-1}(-\sigma)^k$ (with $V(1)=0$); it comes:

$$\begin{matrix}
	W_{i,tors}\simeq\frac{\zd[X]}{(2^l,-1-X)}
\end{matrix}$$

in all cases, of which the order is $2^l$, as expected. The isomorphism described in Case 6 is then proved.
\vx

	\item \underline{ $Char(k)\neq p$:}

\vx

When the characteristic of the residual field $k$ of $K$ is not $p$, since we proved $n=m+b$, we aim to compare $\ft$ with $W_{a,b,n-b,n}$. To do so, we define $\kappa\in\zp^\times$ such that 

$$k_\sigma=1+\kappa p^a$$

and

$$l=\left\{ \begin{matrix}
	 & 2  & & \text{if $p\neq 2$ or $a\neq 1$}\\
	  & & & \\
	 & l \ & \ \text{such that } |\mu_{K(i)}|=2^l & \text{if $p= 2$ and $a= 1$}  
\end{matrix}
\right.$$

(3.7.4) and (3.7.5) then becomes:

$$\pi_F^{1-\sigma}=\xi_F^{p^\delta}=\xi_F^{p^{a-b+l-2}}$$

Then,

$$\left(\pi_F^{p^b}\xi_F^{p^{l-2}\kappa^{-1}}\right)^{1-\sigma}=\xi_F^{p^{a+l-2}}\xi_F^{p^{l-2}\kappa^{-1}(-\kappa p^a)}=1$$

and then we define

$$\pi_K=\pi_F^{p^b}\xi_F^{p^{l-2}\kappa^{-1}}\in\kt$$

which is $G$-invariant and is a uniformizer of $K$. Then, it comes:

$$\frac{\pi_F^{-\sigma+1+\kappa p^a}}{\pi_K^{\kappa p^{a-b}}}=\frac{\xi_F^{p^{a-b+l-2}}\pi_F^{\kappa p^{a}}}{(\pi_F^{p^b}\xi_F^{p^{l-2}\kappa^{-1}})^{\kappa p^{a-b}}}$$

$$\frac{\pi_F^{-\sigma+1+\kappa p^a}}{\pi_K^{\kappa p^{a-b}}}=1$$

This is the fundamental relation of $W_{a,b,m,n}$ for any $m$, replacing $X$ with $\pi_K$, $S$ with $\pi_F$ and $T$ with 1. Then, $\kt$ is a quotient of the formal space

$$W=\frac{\langle X,S\rangle_{\zp[G]}}{\langle S^{-\sigma+1+\kappa p^a}X^{-\kappa p^{a-b}}\rangle_{\zp[G]}}=\frac{M}{D}$$

where $X^\sigma =X$ as usual. The reader may be surprised by the fact that $W$ doesn't depend on $m$, however here $m=n-b$ is indeed the greatest possible value of $m$. In $M$, 

$$\Xi=N_m(S)/X$$

satisfies $$\Xi^{-\sigma+1+\kappa p^a}=q^{N_m}$$

where $q=S^{-\sigma+1+\kappa p^a}X^{-\kappa p^{a-b}}$, and then in $W$, $\Xi$ is $\omega$-isotypical. On the other hand, considering:

$$\omega_n=\left\vert  \frac{\zp[X]}{\langle X^{p^n}-1,-X+1+\kappa p^a\rangle}\right\vert$$

According to (4.3.1) and lemma (3.3.1) applied in $F/K_m$,

$$\omega_n=p^{n-m}\omega_m$$

as soon as $1\le m\le n$. Moreover,

$$\frac{W}{\langle X\rangle}=\frac{\langle S\rangle}{\langle S^{1-\sigma+\kappa p^a}\rangle}\simeq \frac{\zp[X]}{\langle X^{p^n}-1,-X+1+\kappa p^a\rangle}$$

has an order of $\omega_n=p^{n-m}\omega_m$. The corresponding module involving $\ft$ is, if one notes $\phi: W\longrightarrow \ft$ the map sending $S$ to $\pi_F$ and $X$ to $\pi_K$:

$$\frac{\phi(W)}{\langle\phi(X)\rangle}=\frac{\ft}{\langle\pi_K\rangle}$$

It follows:

$$\left\vert\frac{\ft}{\langle\pi_K\rangle}\right\vert=
\left\vert\frac{\ft}{\langle\mu_F\pi_K\rangle}\right\vert\left\vert\frac{\langle\mu_F\pi_K\rangle}{\langle\pi_K\rangle}\right\vert=|Im(\psi)|\times |\mu_F|$$

where $\psi$ is the mapping

$$\psi:\ \begin{pmatrix}
\ft&\longrightarrow &\mathbb{Z}/p^b\mathbb{Z}\\
 w&\mapsto & v(w) \mod p^b
\end{pmatrix}$$

where $v$ is the valuation in $F$ normalized with $v(\pi_F)$=1. Then, 
$$\left\vert Im(\psi)\right\vert\times \left\vert \mu_F\right\vert =p^b\left\vert \mu_F\right\vert  =p^b\omega_m =\omega_m p^{n-m}$$

(see (4.3.1)). Then, $W/\langle X\rangle$ and $\phi(W)/\langle \phi(X)\rangle$ has the same order, that is $\ker(\phi)\subset\langle X\rangle$. Hence, $\ker(\phi)\neq 1$ would imply that $\chi(\ft)=0$ and would be a contradiction, so that finally:

$$W\simeq\ft$$

	\vx
	
\end{enumerate}

\item {\bf Main Theorem}
\vx

Let $F/K$ and $F'/K'$ be 2 cyclic extensions of the same degree $p^n$ of local fields of characteristic zero. Noting $\sigma$ a generator of $\gal(F/K)$ and $\sigma '$ a generator of $\gal(F'/K')$, $\ft$ and $F^{' \times}$ the $p$-completion of the multiplicative groups of $F$ and $F'$, there exists an isomorphism 

$$\Psi: \ft\longrightarrow F^{' \times}$$

such that 

$$\forall w\in \ft,\forall k\in\mathbb{Z},\ \  \Psi(w^{\sigma^k})=(\Psi(w))^{\sigma'^k}$$

if and only if:\vx

\begin{enumerate}
	\item $$\mu_F\simeq\mu_{F'}$$
	
	as $G$-modules, in the following sense:
	
	$$\begin{matrix}
		|\mu_F|=|\mu_{F'}|\\
		\\
		\text{and}\\
		\\
\exists k\in\zp:\ \forall\xi\in\mu_F,\ \ \xi^\sigma=\xi^{k}\ ;\  \forall\xi\in\mu_{F'},\ \ \xi^{\sigma '}=\xi^{k}
	\end{matrix}$$
	
	\item   $$\left\vert\mu_K\cap N_{F/K}(F^{^{\times}})\right\vert=\left\vert\mu_{K'}\cap N_{F'/K'}(F^{'\times})\right\vert$$
	
	\item $$
	\begin{matrix}
	&	\text{  the residual characteristics of $K$ and $K'$ are both different from $p$}\\
	 &	\\
	&	\text{or}\\
	&	\\
	&	\text{the residual characteristics of $K$ and $K'$ are both equal to $p$ and $\dim_{\qp} K=\dim_{\qp} K'$}
	\end{matrix}
	$$
	\vx
	
	Moreover, noting $d=\left\{\begin{matrix}
		& 1 & \text{ if the residual characteristic of $K$ is not $p$}\\
		 & &\\
		 & \dim_{\qp}K & \text{ if the residual characteristic of $K$ equals $p$}
	\end{matrix}\right.$

one has:

	$$\ft\simeq W\oplus\zp[G]^{d-1}$$
	
	\vx
	
	where $W$ is a convenient direct factor of $\ft$ as found in III and discussed in IV.4, with 
	
	$$\chi(W)=\left\{\begin{matrix}
	&		1 & \text{ if the residual characteristic of K is not $p$}\\
	&	 & \\

		&	\chi_{reg}+1 & \text{ if the residual characteristic of K is $p$}
	\end{matrix}\right.$$

and more precisely:

$$W=\frac{\langle X,S,T\rangle_{\zp[G]}}{\langle\mathcal{R}\rangle_{\zp[G]}}=\frac{\zp\oplus\zp[G]\oplus\zp[G]}{\langle\mathcal{R}\rangle}=\frac{M}{D}$$

where $X^\sigma=X$ both in $M$ and $D$ and $\mathcal{R}$ is a fundamental relation in $F/K$, given (with parameters $a,b,m$ defined at the beginning of III and $\kappa$ defined in (3.5.5)) by:

$$\mathcal{R}=\left\{\begin{matrix}
 & S^{-\sigma+1+\kappa p^a}X^{-\kappa p^{a-b}}T^{\sigma^{p^m}-1} & \text{in Cases 1, 2, 3.2, 4, 5}\\
 	 & & \\
 	  & S^{-\sigma+1+\kappa p^a}X^{-\kappa p^{a-b}} & \text{in Case 7, replacing $T$ with 1 in both $M$ and $D$}\\
 	   & & \\
 	   & S^{2^l}X^{-1}& \text{in Case 3.3.a}\\
 	    & & \\
 	    & S^{2^{l-1}(1-\sigma)}T^{-1-\sigma}& \text{in Case 3.3.b}\\
 	     & & \\
 	    & S^{2^l}T^{-1-\sigma}& \text{in Case 6, $-1\in N(\ft)$}\\
 	    & & \\
 	    & S^{2^{l}}X^{-1}T^{-1-\sigma}& \text{in Case 6, $-1\notin N(\ft)$} 
 	 
\end{matrix}
\right.$$

At last, the splitting:

$$W\simeq\frac{W}{\langle X\rangle}\oplus\zp$$

occurs iff $\mu_K\subset N(\ft)$.

\vx

	\vspace{1cm}
	
	Proof:\vx
	
	At first, the condition (a), (b), (c) implies that the parameters $a,b,m$ as defined in the beginning of III are the same for both $F/K$ and $F'/K'$. Note also that they allow to differentiate special case 1 and 2, because in special case 2, $k\equiv -1$ modulo $|\mu_F|$, which is not the case in special case 1 (alternatively, $|\mu_K\cap N(\ft)|$ is not the same in both cases). Then, specifying if $K[\mu_\infty]/K$ is procyclic or not is not necessary - in practice however it is, because it's more easy to know than the precise action of $G$ on $\mu_F$.
	
	Moreover, an isomorphism involving $K$ of which the residual characteristic is $p$, and $K'$ of which the residual characteristic is not $p$, is impossible, because if so:
	
	$$\chi(\ft)=d\chi\{reg\}+1>\chi_{reg}>\chi(F^{'\times})=1$$

	Then, all conditions (a), (b) and (c) are necessarily filled if $\ft\simeq F^{'\times}$.
	Everything else has been proved all along III and IV, except the last assumption. The reader can check each fundamental relation $\mathcal{R}$ on a case-by-case basis (adapting letters denoting variables), however for Case 3.3.a, instead of keeping the original relation 
	
	$$T^{2^l(1-\sigma)}=1$$
	
	with no $X$ both in $M$ and $D$, we noticed that 
	
	$$\frac{\langle T\rangle}{\langle T^{2^l(1-\sigma)}\rangle}\simeq\frac{\langle X,T\rangle}{\langle T^{2^l}X^{-1}\rangle}$$
	
(with $X^\sigma=X$ as usual). The reader should also have a look at the discussion IV.4. The last assumption ($\mu_K\subset N(\ft)$ implies splitting) was mentioned in IV.4.b about $W_{a,b,m,n}$ and then in Cases 1, 2, 3.2, 4, 5, 7. In Cases 3.3.b and Case 6 (when -1 is a norm), $X$ just doesn't appear in $\mathcal{R}$ and the splitting is then obvious.

About the converse (splitting implies $\mu_K\subset N(\ft)$), we observe that Cases 1, 2, 3.2, 3.3.b are not concerned ($\mu_K\subset N(\ft)$), in Cases 4, 5, one has modulo X: $N(S)=\xi_K$ an element of order $p^a$ which is wrong in $\ft$ if $\mu_K\not\subset N(\ft)$. In Case 3.3.a, $X=S^{2^l}$ with $l\ge 2$ cannot generate a direct factor in $\ft$. In Case 6 when $-1$ is not a norm, $\mathcal{R}$ modulo $X$ gives the fundamental relation obtained when $-1$ is a norm and a splitting then implies the existence of an invariant element of order 2 being a norm, which is not the case since $-1$ is not a norm. Finally in Case 7, IV.4.g states that $|W/\langle X\rangle|=\omega_n$ while $|\mu_F|=|\ft_{tors}|=\omega_m$, hence a splitting implies $m=n$ that is $\mu_K\subset N(\ft)$. All this proves the converse. Note however that this is valid only for these $X$ as chosen all along this paper: in Case 7, when $F/K$ is totally ramified, one has $\mu_K=\mu_F$ and $X=\xi_K$ generates $\kt/N(\ft)$, $\langle X\rangle$ being a direct factor in $\ft$.

\end{enumerate}

	\end{enumerate}

\vspace{1cm}

{\bf \large V Conclusion } \vspace{1cm}

The main motivation of this paper was to determine whether basic parameters were enough or not to describe completely the action of $G$ in a cyclic extension $F/K$ of local fields of characteristic zero. So, the answer is yes. Case 5 was hard to perform, but turned out to be the most important case and underlines the role of the character $\omega$. Moreover, it allowed to find the parameterized space $W_{a,b,m,n}$, which was completely unexpected and applies in almost all situations. Another important point is that $W$ - and $W/\langle X\rangle$ when splitting is possible - concentrates all the interesting things about $F/K$. The author is convinced that it has a lot to offer.

\vspace{2cm}

\makeatletter
\renewcommand{\@biblabel}[1]{}
\makeatother

\end{document}